\documentclass[10pt]{amsart}
\usepackage{amsfonts,psfig,pb-diagram}

\newcommand{\bZ}{\mbox{${\mathbb Z}$}}

\newcommand{\bC}{\mbox{${\mathbb C}$}}

\newcommand{\setp}[2]{\mbox{$\{#1\::\:#2\}$}}

\newcommand{\case}[5]{#1#2 \left\{ \begin{array}{ll} #3 &\mbox{if}\spa #4 \\ #5 &\mbox{otherwise\,.} \end{array} \right.}

\newcommand{\spa}{\;\:}

\newcommand{\fs}{{\mathfrak S}}

\newtheorem{prop}[equation]{Proposition}

\newtheorem{alg}[equation]{Algorithm}
\newtheorem{thm}[equation]{Theorem}
\newtheorem{rem}[equation]{Remarks}
\newtheorem{cor}[equation]{Corollary}
\newtheorem{lem}[equation]{Lemma}
\newtheorem{exa}[equation]{Example}

\newtheorem{rema}[equation]{Remark}

\makeatletter
\let\choose\@@choose
\let\cal\mathcal
\makeatother

\makeatletter\@addtoreset{equation}{section}\makeatother

\makeatletter
\let\cal\mathcal
\makeatother

\begin{document}
\bibliographystyle{plain}

\title[Combinatorial Formulas for Schubert polynomials]{A Unified Approach to Combinatorial Formulas for Schubert Polynomials}
\author{Cristian Lenart}
\address{Department of Mathematics and Statistics, SUNY at Albany, Albany, NY 12222}
\email{lenart@csc.albany.edu}
\thanks{Most of this work was done while the author was supported by Max-Planck-Institut f\"ur Mathematik.}

\begin{abstract}
Schubert polynomials were introduced in the context of the geometry of flag varieties. This paper investigates some of the connections not yet understood between several combinatorial structures for the construction of Schubert polynomials; we also present simplifications in some of the existing approaches to this area. We designate certain line diagrams known as rc-graphs as the main structure. The other structures in the literature we study include: semistandard Young tableaux, Kohnert diagrams, and balanced labelings of the diagram of a permutation. The main tools in our investigation are certain operations on rc-graphs, which correspond to the coplactic operations on tableaux, and thus define a crystal graph structure on rc-graphs; a new definition of these operations is presented. One application of these operations is a straightforward, purely combinatorial proof of a recent formula (due to Buch, Kresch, Tamvakis, and Yong), which expresses Schubert polynomials in terms of products of Schur polynomials. In spite of the fact that it refers to many objects and results related to them, the paper is mostly self-contained.
\end{abstract}

\maketitle\vspace{-1cm}


\section{Introduction}

Schubert polynomials are relevant to the geometry of flag varieties, and their study leads to some subtle combinatorics, generalizing the combinatorics of semistandard Young tableaux (SSYT) related to Schur polynomials (which are special cases). Schubert polynomials were defined and investigated by Lascoux and Sch\"{u}tzenberger in a series of papers, including \cite{lasdcg,lasps,lassha,lasfps,lastns,lasksb}. The combinatorics related to Schubert polynomials attracted the interest of many people in algebraic combinatorics (see the references). Several combinatorial constructions of Schubert polynomials appeared in the literature in the last decade, together with many generalizations of them (such as generalizations to other Lie groups, to the $K$-theory and the quantum cohomology of flag varieties etc.). While all these generalizations appeared, the connections between the existing structures for the construction of Schubert polynomials themselves were not fully understood. In this paper we attempt to elucidate some of the connections not yet understood, and to present simplifications in some of the existing approaches to this area. 

We start with a brief introduction to the theory of Schubert polynomials; for more information, we refer the reader to \cite{fulyt,macnsp,manfsp}. Let $Fl_n$ be the variety of complete flags $(0=V_0\subset V_1\subset\ldots \subset V_{n-1}\subset V_n=\bC^n)$ in $\bC^n$. This is an irreducible algebraic variety of complex dimension $\binom{n}{2}$.  Its integral cohomology ring $H^*(Fl_n)$ is isomorphic to $\bZ[x_1,\ldots,x_n]/I_n$, where $I_n$ is the ideal generated by symmetric polynomials in $x_1,\ldots,x_n$ with constant term 0; the elements $x_i$ are identified here with the Chern classes of the dual line bundles $L_i^*$, where $L_i:=V_i/V_{i-1}$ are tautological line bundles. The variety $Fl_n$ is a disjoint union of cells indexed by permutations $w$ in the symmetric group $S_n$. Their closures are the {\em Schubert varieties} $X_w$, of complex dimension $l(w)$; here $l(w)$ denotes the length of $w$, that is, the number of its inversions. The Schubert polynomial $\fs_w(X)$ (defined below) is a certain polynomial representative for the cohomology class corresponding to $X_w$. It is a homogeneous polynomial in $x_1,\ldots,x_{n-1}$ of degree $l(w)$ with nonnegative integer coefficients. 

Note that the above structure theorem for $H^*(Fl_n)$ holds, more generally, for flag bundles $Fl(V)$ corresponding to a vector bundle $V$ of finite rank $n$ on a variety $X$. In this case, $H^*(Fl(V))$ is a free module over $H^*(X)$ with basis given by Schubert polynomials. 

The construction of Schubert polynomials is based on {\em divided difference operators} $\partial_i$, which were originally introduced in the context of generalized flag varieties \cite{bggscc,kaknhr}. As operators on $\bZ[x_1,x_2,\ldots]$, these are defined as follows:
\begin{equation}\label{defpi}
\partial_i:=\frac{1-s_i}{x_i-x_{i+1}}\,;
\end{equation}
here $s_i$ is the adjacent transposition $t_{i,i+1}$ acting in the obvious way on the polynomial algebra, $1$ is the identity on this algebra, and $x_i$ denotes multiplication by this variable. 

Schubert polynomials are defined inductively for each $w$ in $S_n$ by setting $\fs_{w_0^n}(X):=x_1^{n-1}x_2^{n-2}\ldots x_{n-1}$ for $w_0^n:=n(n-1)\ldots 1$ the longest permutation in $S_n$ (we use the one-line notation for permutations throughout, unless otherwise specified), and by letting
\begin{equation}\label{rec1}
\fs_{ws_i}(X):=\partial_i(\fs_w(X))\,,\;\;\;\;\;\mbox{if}\spa l(ws_i)=l(w)-1\,.
\end{equation}
The definition of Schubert polynomials does not depend on the chosen chain in the weak order on $S_n$ from $w_0^n$ to $w$ because the operators $\partial_i$ satisfy the braid relations
\begin{align*}
&\partial_i\partial_j=\partial_j\partial_i\;\;\;\;\mbox{if}\spa |i-j|\ge 2\,,\\
&\partial_i\partial_{i+1}\partial_i=\partial_{i+1}\partial_i\partial_{i+1}\,.
\end{align*}
In addition, we have $\partial_i^2=0$. Also note that Schubert polynomials have a stability property, in the sense that their definition does not depend on $n$; this allows one to define them for $w$ in $S_{\infty}:=\bigcup_n S_n$, where $S_n\hookrightarrow S_{n+1}$ is the usual inclusion. In fact, the Schubert polynomials $\fs_w(X)$ for $w$ in $S_{\infty}$ form a basis of $\bZ[x_1,x_2,\ldots]$.

As we mentioned above, there are several combinatorial constructions of Schubert polynomials. We designate certain line diagrams for permutations introduced in \cite{fakybe} and called {\em rc-graphs} in \cite{babrcg} (see Section 2) as the main combinatorial structure for the construction of Schubert polynomials. The other structures in the literature we study include: semistandard Young tableaux \cite{lasdcg,lastns, lasksb,raskpf} (Sections 3-5), Kohnert diagrams \cite{kohwpt, windrg, winspk} (Section 6), and balanced labelings of the diagram of a permutation \cite{fgrbls} (Section 7). We do not mention here the construction of Schubert polynomials in terms of certain increasing labeled chains in Bruhat order on the symmetric group \cite{basssf}; the connection between these chains and rc-graphs is explored in \cite{lasssp}. The main tools in our investigation are certain operations on rc-graphs, which correspond to the coplactic operations on tableaux, and thus define a crystal graph structure on rc-graphs. A new definition of these operations is presented, and their properties are studied in Section 3. Among the several applications of these operations is a straightforward, purely combinatorial proof (in Section 5) of a recent formula \cite{bktspq}, which expresses Schubert polynomials in terms of products of Schur polynomials.


{\bf Acknowledgements.} The author is grateful to Sara Billey, Mark Shimozono, Ezra Miller, Anders Buch, Harry Tamvakis, and Alexander Yong for helpful discussions, as well as to one of the referees for providing Remark \ref{referee}. 

\section{RC-Graphs}

RC-graphs are combinatorial objects associated to permutations $w=w_1\ldots w_n$ in the symmetric group $S_n$. We let $n$ and $w$ be fixed throughout, unless otherwise specified. RC-graphs can be defined as certain subsets of
\[\setp{(i,j)\in\bZ_{>0}\times\bZ_{>0}}{i+j\le n}\,.\]
Given an arbitrary such subset of pairs $R$, we can define a linear order on it by setting
\begin{equation}\label{ordcross}(i,j)\le (i',j') \;\;\Longleftrightarrow\;\;  (i< i') \spa\mbox{or} \spa (i=i'\spa\mbox{and}\spa j\ge j')\,.\end{equation}
Let $(i_k,j_k)$ be the $k$th pair in this linear order, and let $a_k:=i_k+j_k-1$. Consider the sequences (words)
\[{\rm red}(R):=a_1a_2\ldots\,,\;\;\;\;\; {\rm comp}(R):=i_1i_2\ldots\,.\]
We call $R$ an rc-graph associated to $w$ if ${\rm red}(R)$ is a reduced word for $w$; in this case, ${\rm comp}(R)$ is called a {\em compatible sequence} with ${\rm red}(R)$. We denote the collection of rc-graphs associated to a permutation $w$ by ${\mathcal R}(w)$, and the permutation corresponding to a given rc-graph $R$ by $w(R)$. 

An alternative way to define ${\mathcal R}(w)$ is to consider all reduced words $a_1\ldots a_{l(w)}$ for $w$, and for each such word to consider all weakly increasing sequences $i_1\ldots i_{l(w)}$ of positive integers satisfying the ``compatibility'' conditions:
\begin{align}\label{comp1}
&a_k<a_{k+1} \mbox{ implies } i_k<i_{k+1} \mbox{ for all } 1\le k<l(w)\,;\\
&i_k\le a_k \mbox{ for all } 1\le k\le l(w)\,. \label{comp2}
\end{align}
Note that some reduced words may have no compatible sequences associated to them.

We can represent an rc-graph $R$ graphically as a planar history of the inversions of $w$ (see Example \ref{exrc}). To this end, we draw $n$ lines going up and to the right, such that the $i$th line starts at position $(i,1)$ and ends at position $(1,w_i)$; here the positions are numbered as in a matrix, whence the first index is the row index, and the second one is the column index. The rule for constructing the line diagram is the following: two lines meeting at position $(i,j)$ cross if  $(i,j)$ is in $R$; otherwise, they avoid each other at that position. Throughout this paper, we will use interchangeably the three representations of an rc-graph mentioned above (that is, as a subset of $\bZ_{>0}\times\bZ_{>0}$, as a pair of a reduced word and a compatible sequence, and as a line diagram).

Given an rc-graph $R$, we define
\[x^R:=\prod_{(i,j)\in R} x_i\,.\]
According to \cite{babrcg,bjsscp,fakybe,fasspn,lassha}, the Schubert polynomial $\fs_w(X)$ indexed by $w$ can be expressed as
\begin{equation}\label{schub1}
\fs_w(X)=\sum_{R\in{\mathcal R}(w)}x^R\,.
\end{equation}
Note that this formula was recently recovered in \cite{kamggs} in the algebraic and geometric context of determinantal ideals specified by rank conditions, based on the combinatorics of initial ideals for certain natural term orders. 

Bergeron and Billey \cite{babrcg} defined certain moves on rc-graphs, which they called {\em ladder} and {\em chute} moves. These are the moves $L_{ij}$ and $C_{ij}$ of the following type:
\[
\begin{array}{c}
\mbox{\psfig{file=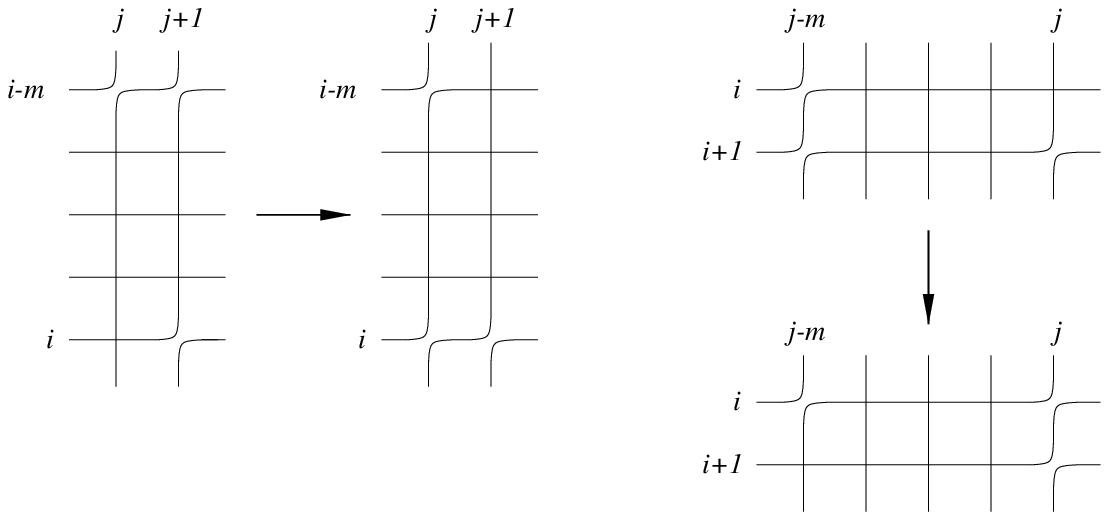}}
\end{array}
\]
Note that $L_{ij}$ moves the cross in position $(i,j)$ to some position $(i-m,j+1)$; we will denote it by $(i,j)\rightarrow(i-m,j+1)$.

Recall that a permutation $w$ is determined by its {\em code}, denoted ${\rm code}(w)$, which is the sequence $(c_1,\ldots,c_{n-1})$ defined by $c_i=|\setp{j}{j>i\;\;\mbox{and}\;\; w_j<w_i}|$. Let us also recall the rc-graphs $R_{\rm bot}(w)$ and $R_{\rm top}(w)$ defined by
\[R_{\rm bot}(w):=\setp{(i,j)}{j\le c_i(w)}\,,\;\;\;\;\;\;\;\; R_{\rm top}(w):=\setp{(i,j)}{i\le c_j(w^{-1})}\,.\]
Bergeron and Billey showed that ${\mathcal R}(w)$ is the transitive closure of $\{R_{\rm bot}(w)\}$ under ladder moves $L_{ij}$ which move the rightmost cross in a row (i.e., $(i,l)$ is not in the rc-graph to which $L_{ij}$ is applied if $l>j$); we call these moves, simply, {\em L-moves}. There is a dual property, namely that ${\mathcal R}(w)$ is the transitive closure of $\{R_{\rm top}(w)\}$ under chute moves $C_{ij}$ which move the bottom cross in a row (i.e., $(k,j)$ is not in the rc-graph to which $C_{ij}$ is applied if $k>i$). 

\begin{exa}\label{exrc}{\rm Here are two rc-graphs associated to the permutation $w=215463$; the first is just $R_{\rm bot}(w)$. 
\[\begin{array}{c}
\mbox{\psfig{file=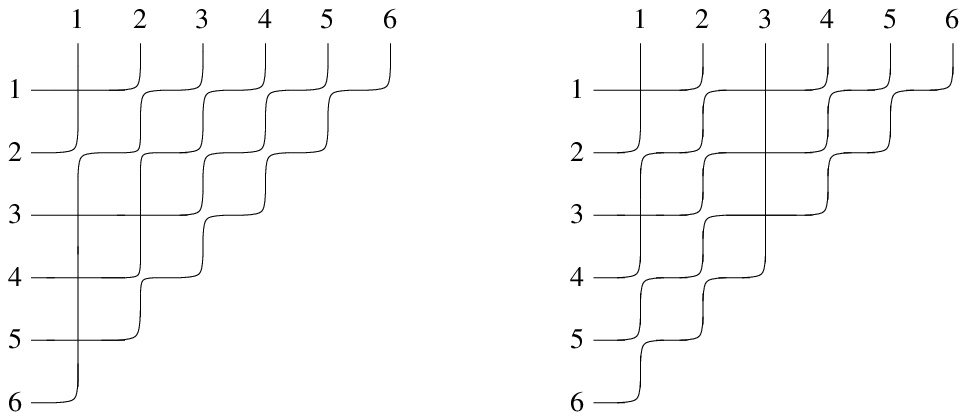}}
\end{array}
\]
}\end{exa}

It is useful to have a standard way of constructing a given rc-graph $R$ for $w$ from $R_{\rm bot}(w)$ via L-moves. This can be done by applying the algorithm below, based on the Proposition following it.

\begin{alg}\label{std}$\;\;$\newline
\noindent set $R':=R_{\rm bot}(w)$;\newline
for $i$ from 1 to $n-1$ do\newline
 \indent for $j$ from $n-i$ downto 1 do\newline
\indent\indent if $(i,j)\in R\setminus R'$ then\newline
\indent\indent\indent find the position $(k,l)$ where the lines of $R'$ avoiding each other at $(i,j)$ cross;\newline
\indent\indent\indent set $R':=R'\cup\{(i,j)\}\setminus\{(k,l)\}$;\newline
\indent\indent fi;\newline
\indent od;\newline
od.
\end{alg}

\begin{prop} In the above algorithm, each cross in position $(k,l)$ is moved into position $(i,j)$ by a sequence of L-moves.
\end{prop}
\begin{proof}
Assume that the assertion is true for all previous moves. This implies that all crosses are left justified in the part of the rc-graph consisting of positions $(p,q)>(i,j)$ in the order (\ref{ordcross}), that is, the unshaded region in the figure below. 
\[
\begin{array}{c}
\mbox{\psfig{file=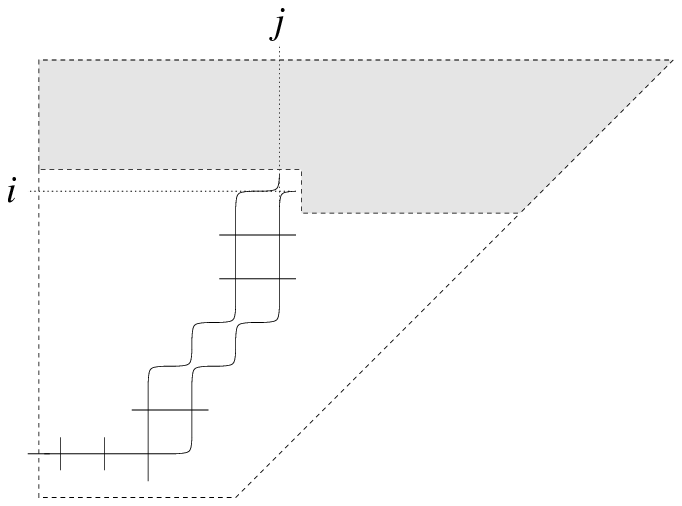}}
\end{array}
\]
It is not hard to check, based on this observation, that the two lines avoiding each other at position $(i,j)$ must follow a route of the form indicated in the figure below. Furthermore, there are no crosses to the right of the positions $(p,q)>(i,j)$ where the two lines avoid each other (if any). 
\end{proof}

\begin{exa}\label{exmv}{\rm
The standard sequence of L-moves for the second rc-graph in Example \ref{exrc} is 
$(4,1)\rightarrow(2,2)$, $(2,2)\rightarrow(1,3)$, $(3,2)\rightarrow(2,3)$, $(5,1)\rightarrow(4,2)$, $(4,2)\rightarrow(3,3)$.
}\end{exa}

The main feature of this construction is that the sequence of ``final'' positions of the crosses which are moved is increasing in the order defined in (\ref{ordcross}). We can use this idea to improve the algorithm of Bergeron and Billey based on L-moves, by constructing ${\mathcal R}(w)$ such that each rc-graph is obtained precisely once. This procedure can be illustrated by constructing a tree whose root corresponds to $R_{\rm bot}(w)$, and whose edges correspond to applying L-moves. For each rc-graph $R$ (corresponding to a given vertex in the tree), we use the variable ${\rm last}(R)$ to record the position of the last cross moved to a ``final'' position. We also use the variable ${\rm current}(R)$ to record the position of the last moved cross, possibly to a ``temporary'' position. This cross might or might not be moved by a subsequent L-move; if it is not moved, then it becomes the last cross moved to a ``final'' position. For $R_{\rm bot}(w)$, these two variables are undefined. We can now describe the branching rule at an rc-graph $R$ in our tree; the order on pairs used is the one defined in (\ref{ordcross}).

\begin{alg}\label{allrc}$\;\;$\newline
\noindent for all L-moves $P\rightarrow Q$ in $R$ do\newline
\indent if $P={\rm current}(R)$ then\newline
\indent\indent if ${\rm last}(R)$ is undefined or $Q>{\rm last}(R)$ then\newline
\indent\indent\indent construct $R'$ by applying the move $P\rightarrow Q$ to $R$;\newline
\indent\indent\indent let ${\rm current}(R'):=Q$, ${\rm last}(R'):={\rm last}(R)$;\newline
\indent\indent fi;\newline
\indent else \newline
\indent\indent if ${\rm current}(R)$ is defined then\newline
\indent\indent\indent if $Q>{\rm current}(R)$ then\newline
\indent\indent\indent\indent construct $R'$ by applying the move $P\rightarrow Q$ to $R$;\newline
\indent\indent\indent\indent let ${\rm current}(R'):=Q$, ${\rm last}(R'):={\rm current}(R)$;\newline
\indent\indent\indent fi;\newline
\indent\indent else \newline
\indent\indent\indent construct $R'$ by applying the move $P\rightarrow Q$ to $R$;\newline
\indent\indent\indent let ${\rm current}(R'):=Q$;\newline
\indent\indent fi;\newline
\indent fi;\newline
od.
\end{alg}

We believe this algorithm to be the most efficient one for constructing ${\mathcal R}(w)$ for $w$ in $S_n$. Its complexity is clearly $O(nc)$, where $c$ is the cardinality of ${\mathcal R}(w)$. Other algorithms for constructing ${\mathcal R}(w)$, and thus the multiset of monomials in $\fs_w(X)$, appeared in \cite{lasssp} and \cite{milmrc}. As discussed in \cite{lasssp}, their complexity is $O(nlc)$, where $l:=l(w_0^n)-l(w)=\binom{n}{2}-l(w)$.

\section{Double Crystal Graphs}

In this section we study in more detail the double crystal graphs introduced in \cite{lasdcg}, with an eye towards applications to the theory of Schubert polynomials; the results are then used in Sections 4-6. More precisely, we study the main properties of the generalizations of the coplactic operations on words to biwords in the cases of both the {\em plactic} and the {\em nilplactic monoid}. We present a simple way of deducing these properties from those of the classical coplactic operations. The nilplactic case is then easily reduced to the plactic one via the {\em plactification map} defined in \cite{raspla}. The main feature of our approach, as opposed to the one in \cite{lasdcg}, is that it is centered around the concept of rc-graph (for instance, we define coplactic operations directly on rc-graphs via $r$-pairing of crosses in two consecutive rows). Indeed, the construction of Schubert polynomials in terms of rc-graphs (\ref{schub1}) is mentioned as a corollary in the above paper, whereas, for us, it is the starting point. 

We start by recalling the setup, and we refer to \cite{eagbt,fulyt,lasdcg,lltpm} for more details, including the standard facts about the combinatorics of Young tableaux. This combinatorics is underlied by the plactic monoid, which is the monoid of words on a given alphabet modulo the {\em Knuth} (or {\em plactic}) {\em equivalence}; this, in turn, is the transitive closure of the relations
\begin{align}\label{knuth}
\ldots ikj\ldots\,&\sim_{K}\,\ldots kij\ldots\spa\spa\mbox{for $i\le j<k$}\\
\ldots jik\ldots\,&\sim_{K}\,\ldots jki\ldots\spa\spa\mbox{for $i< j\le k$}\,.\nonumber
\end{align}
Let $X$ and $Y$ be two totally ordered alphabets; for convenience, we take $X=Y:=\bZ_{>0}$. We consider {\em biwords}, that is, two-rowed arrays
\[w:=\binom{u}{v}=\binom{u_1\ldots u_k}{v_1\ldots v_k}\]
with $u_i$ in $X$ and $v_j$ in $Y$. Recall that $w$ is said to be in {\em lexicographic order} if its biletters $\binom{u_i}{v_i}$ satisfy
\[(u_i<u_{i+1})\spa\mbox{or}\spa(u_i=u_{i+1}\spa\mbox{and}\spa v_i\le v_{i+1})\spa\mbox{for}\spa i=1,\ldots,k-1\,.\]
On the other hand, $w$ is said to be in {\em antilexicographic order} if its biletters satisfy 
\[(u_i>u_{i+1})\spa\mbox{or}\spa(u_i=u_{i+1}\spa\mbox{and}\spa v_i\le v_{i+1})\spa\mbox{for}\spa i=1,\ldots,k-1\,.\]
Given a biword $w=\binom{u}{v}$, we let $w^{\rm rev}:=\binom{u^{\rm rev}}{v^{\rm rev}}$ be the biword obtained by considering the corresponding reverse words. Furthermore, we let $w'=\binom{v'}{u'}$ be the unique biword obtained by rearranging the biletters of $\binom{v}{u}$ such that $(w')^{\rm rev}$ is in lexicographic order. From now, we will only consider words with no two identical biletters, so the corresponding lexicographic and antilexicographic orders are strict. 

Given a biword $w=\binom{u}{v}$ with $w^{\rm rev}$ in antilexicographic order, we can associate with it a pair of SSYT of conjugate shapes $(P(w),Q(w))$; here $P(w)$ is obtained by applying Schensted insertion to the word $v$, while $Q(w)$ is obtained by {\em conjugate placing} the corresponding letters of $u$. This is a variation of the Robinson-Schensted-Knuth (RSK) correspondence, which sets up a bijection between matrices of zeros and ones of finite support and pairs of SSYT of conjugate shapes. On the other hand, we can associate a similar pair of SSYT to $w'$, denoted by $(P(w'),Q(w'))$, with the only exception that $Q(w')$ is obtained by {\em conjugate sliding} the corresponding entries of $v'$. According to the symmetry property of the mentioned variation of the RSK correspondence, we have that
\begin{equation}\label{symm}P(w)=Q(w')\spa\mbox{and}\spa P(w')=Q(w)\,.\end{equation}
For more information on the above constructions we refer to \cite{fulyt} Appendix A.4.3. 

The concepts of left and right strings for biwords in lexicographic order were defined in \cite{lasdcg}; hence, we have the corresponding generalizations to biwords of the coplactic operations on words, which we denote by $e_r^-$, $f_r^-$, and $e_r^+$, $f_r^+$. We now define these operations on biwords $w=\binom{u}{v}$ with $w^{\rm rev}$ in antilexicographic order. Write the word $u$ in the form $1^{l_1}\ldots p^{l_p}$, where $i^l$ is the word with $l$ letters equal to $i$, and $v$ in the form ${\rm col}_1(w)\ldots {\rm col}_p(w)$, where ${\rm col}_i(w)$ are strictly decreasing sequences (columns) of lengths $l_i$. Given a pair of columns $(u,v)$, we can apply jeu de taquin or inverse jeu de taquin to produce two different pairs of columns, which are denoted by ${\rm jdt}(u,v)$ and ${\rm jdt}^{-1}(u,v)$, respectively (see below). With this in mind, we define
\[e_r^+(w):=\binom{u\!\uparrow}{v\!\uparrow}\,,\spa f_r^+(w):=\binom{u\!\downarrow}{v\!\downarrow}\,,\]
where 
\[u\!\uparrow:=1^{l_1}\ldots r^{l_r+1}(r+1)^{l_{r+1}-1}\ldots p^{l_p}\,,\spa u\!\downarrow:=1^{l_1}\ldots r^{l_r-1}(r+1)^{l_{r+1}+1}\ldots p^{l_p}\,,\]
and 
\begin{align*}&v\!\uparrow:={\rm col}_1(w)\ldots{\rm jdt}({\rm col}_r(w),{\rm col}_{r+1}(w))\ldots{\rm col}_p(w)\,,\\ &v\!\downarrow:={\rm col}_1(w)\ldots{\rm jdt}^{-1}({\rm col}_r(w),{\rm col}_{r+1}(w))\ldots{\rm col}_p(w)\,.\end{align*}

Let us now define precisely what we mean by jeu de taquin and inverse jeu de taquin on two columns
\[u:=(u_s>u_{s-1}>\ldots >u_1)\spa\mbox{and}\spa v:=(v_t>v_{t-1}>\ldots >v_1)\,.\]
We place the two columns side by side, in a skew tableau denoted $T(u,v)$, such that there is a maximum overlap between them. In other words, we find the unique $i$ with ${\rm max}(0,t-s)\le i\le t$ such that the following two conditions are satisfied (see the figure below):
\begin{enumerate}
\item $u_j\le v_{i+j}$ for $1\le j\le t-i$;
\item if $i>\max(0,t-s)$, then $u_k>v_{i+k-1}$ for some $k$ with $1\le k\le t-i+1$.
\end{enumerate}
If $i=0$, then we define ${\rm jdt}(u,v):=\emptyset$, and if $i=t-s$, then we define ${\rm jdt}^{-1}(u,v):=\emptyset$. If such a situation occurs, the result of applying $e_r^+$ or $f_r^+$ is defined to be the empty biword $\emptyset$. 

\begin{exa}{\rm Here is an example of jeu de taquin on two columns.
\[
\begin{array}{c}
\mbox{\psfig{file=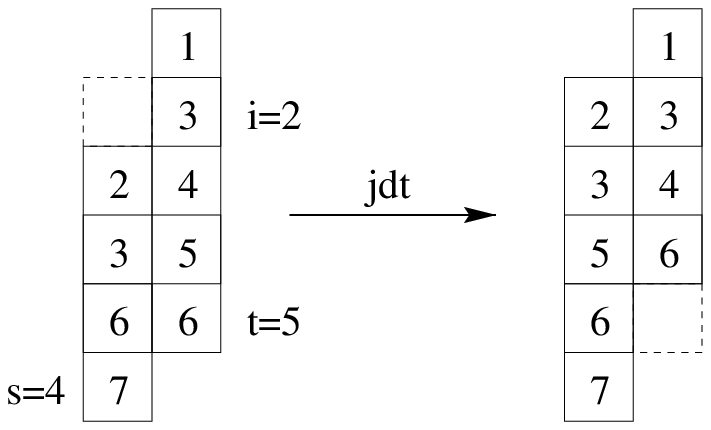}}
\end{array}
\]
}\end{exa}

The following result is implicitly used in \cite{lasdcg}.

\begin{lem}\label{uniqueright}
There is at most one pair of columns of lengths $s+1$ and $t-1$ in the plactic class of the word $uv$ (the concatenation of $u$ and $v$). If such a pair exists, then it is precisely ${\rm jdt}(u,v)$. Similarly for ${\rm jdt}^{-1}(u,v)$. 
\end{lem}
\begin{proof} It is well-known (see e.g. Corollary 2 on p. 62 in \cite{fulyt}) that the number of skew tableaux of a fixed shape in the plactic class of a given SSYT $T^*$ depends only on the shape of that tableau. Furthermore, this number is a Littlewood-Richardson coefficient; thus, we get the same number if we transpose all tableaux involved. Now assume that ${\rm jdt}(u,v)\ne\emptyset$. By jeu de taquin, $T(u,v)$ is in the plactic class of a tableau of shape $(s+i,t-i)^T$ with $i>0$. Now let $T^*:=U(s+i,t-i)$, where $U(\lambda)$ denotes the tableau of shape $\lambda$ with all entries in row $i$ equal to $i$. It is easy to check that there is a unique skew tableaux of shape $(s+t,t-1)/(t-1)$ in the plactic class of $T^*$; indeed, all the entries in its first rows must be equal to 1. The case ${\rm jdt}(u,v)=\emptyset$ is obvious. Hence the Lemma follows. A similar reasoning works for the second part of the Lemma. \end{proof}

Let us now recall the definition of $e_r^-$ and $f_r^-$ on biwords. In fact, this works for any biword $w$, so no assumption on the order of the biletters of $w$ is needed. We simply set
\[e_r^-(w):=\binom{u}{e_r(v)}\,,\spa f_r^-(w):=\binom{u}{f_r(v)}\,,\]
where $e_r$ and $f_r$ are the usual coplactic operations on words, whose definition is based on the concept of {\em r-pairing}. If $e_r(v)=\emptyset$ (resp. $f_r(v)=\emptyset$), then we set $e_r^-(w):=\emptyset$ (resp. $f_r^-(w):=\emptyset$).

Recall that the $r$-pairing of a word $u$ is just a set of indexed pairs (called $r$-pairs) $(u_i,u_j)$ such that $i<j$, $u_i=r+1$, and $u_j=r$. View each $r+1$ (resp. $r$) as a left (resp. right) paranthesis, ignoring the other letters of $u$. The $r$-pairs of $u$ are precisely the matched parantheses. The subword of unpaired $r$'s and $r+1$'s is necessarily of the form $r^s(r+1)^t$. Define $e_r(u)$ to be the word obtained by changing the leftmost unpaired $r+1$ to $r$ if $t>0$; otherwise, define it to be the empty word $\emptyset$. Similarly, $f_r(u)$ is defined by changing the rightmost unpaired $r$ to $r+1$ if $s>0$. The operations $f_r$ define a directed graph $\Gamma$ on words in $X^*$, whose edges are the ordered pairs $(u,f_r(u))$, for $f_r(u)\ne\emptyset$; this is known as a {\em crystal graph}. 

\begin{rema}\label{uniqueleft}{\rm Consider an arbitrary word $u$, and denote by $\overline{u}$ its subword consisting of letters $r,r+1$. It is easy to check that there is at most one word $u\!\uparrow$ obtained from $u$ by changing precisely one letter $r+1$ to $r$ such that $Q(\overline{u\!\uparrow})=Q(\overline{u})$. If such a word exists, then it is precisely $e_r(u)$. Similarly for $f_r(u)$.}
\end{rema}

Let us also recall at this point the operation $\sigma_r$, where $\sigma_r(u)$ is defined by replacing the subword $r^s(r+1)^t$ with $r^t(r+1)^s$. This gives an action of the symmetric group on words. The main properties of the operations $e_r$, $f_r$, $\sigma_r$ are summarized below.

\begin{thm}\label{erfrsr1}\cite{lasmp,lltpm} Let $h$ denote anyone of the operations $e_r$, $f_r$, $\sigma_r$.\\
\indent {\rm (1)} Let $u$ the the row (resp. column) word of a skew tableau $T$. Then $h(u)$ is the row (resp. column) word of a skew tableau of the same shape; this will be denoted by $h(T)$. \\
\indent {\rm (2)} Let $u$ and $v$ be congruent words in the plactic monoid. Then $h(u)$ is congruent to $h(v)$. In particular, $P(h(u))=h(P(u))$. \\
\indent {\rm (3)} Let $u$ be a word, and assume that $h(u)\ne \emptyset$. Then $Q(h(u))=Q(u)$.\\
\indent {\rm (4)} The connected components of $\Gamma$ are the coplactic classes in $X^*$. Two coplactic classes are isomorphic as subgraphs of $\Gamma$ if and only if they are indexed by standard tableaux of the same shape.
\end{thm}

\begin{rem}{\rm {\rm (1)} The first property implies that if $w$ or $w^{\rm rev}$ are in lexicographic or antilexicographic order, the same is true for $e_r^-(w)$ and $f_r^-(w)$ (assuming they are nonempty).\\ 
\indent {\rm (2)} Clearly, we have 
\[e_r^+(w)=\widetilde{w}\spa\mbox{if and only if}\spa f_r^+(\widetilde{w})=w\,,\]
and similarly for $e_r^-$, $f_r^-$. 
}\end{rem} 

We now present the relationship between the operations $e_r^-$, $f_r^-$ on the one hand, and $e_r^+$, $f_r^+$, on the other hand. This is mentioned in \cite{lasdcg} in a slightly different context (of row words, rather than column words). The author argues that the symmetry in the two alphabets $X$ and $Y$ allows one to exchange the left and right strings. We make this argument explicit in the case of column words.

\begin{prop}\label{plusminus} For any biword $w$ with $w^{\rm rev}$ in antilexicographic order, we have
\[e_r^-(w')=(e_r^+(w))'\,,\]
and similarly for $f_r^-$, $f_r^+$.
\end{prop}
\begin{proof}
Consider a biword $w=\binom{u}{v}$ and its subword, denoted $\overline{w}$, consisting of biletters with top component equal to $r$ or $r+1$. Denote $\overline{w}:=\binom{\overline{u}}{\overline{v}}$ and $\overline{w}':=\binom{\overline{v}'}{\overline{u}'}$. Given a standard tableau $T$ with $n$ entries and a word $v=v_1\ldots v_n$, we denote by ${\rm evac}(T)$ the {\em evacuation} of $T$, and by ${\rm subs}(v,T)$ the array obtained by replacing each entry $i$ in $T$ by $v_i$. 

Assume that $e_r^-(w')\ne\emptyset$, and let $\widetilde{w}$ be the biword with $\widetilde{w}^{\rm rev}$ in antilexicographic order satisfying $\widetilde{w}'=e_r^-(\overline{w}')$. According to (\ref{symm}) and Theorem \ref{erfrsr1} (3), we have
\begin{align*}
P(\widetilde{w})&=Q(e_r^-(\overline{w}'))={\rm subs}(\overline{v}',{\rm evac}(Q^T(e_r(\overline{u}'))))\\
&={\rm subs}(\overline{v}',{\rm evac}(Q^T(\overline{u}')))=Q(\overline{w}')=P(\overline{w})\,.
\end{align*}
By Lemma \ref{uniqueright}, it follows that $\widetilde{w}=e_r^+(\overline{w})$, and thus the desired identity is proved. The only thing we still need to prove is that $e_r^+(w)\ne\emptyset$ implies $e_r^-(w')\ne\emptyset$. This follows using similar ideas, except that we need Remark \ref{uniqueleft} instead of Lemma \ref{uniqueright}. 
\end{proof}

It turns out that it is convenient to identify in the obvious way a biword $w=\binom{u}{v}$ for which $w^{\rm rev}$ is in antilexicographic order with a finite subset $D$ of $\bZ_{>0}\times\bZ_{>0}$ (each biletter $\binom{u_i}{v_i}$ corresponds to a pair $(u_i,v_i)$); this is called a {\em diagram}, and its elements are called {\em cells}. We also write $u={\rm row}(D)$ and $v={\rm col}(D)$. Note that biwords are usually identified with (0,1)-matrices of finite support, which is essentially the same thing, although our approach is slightly more convenient. On the other hand, the matrix representation is sometimes useful, so we will use all the above representations interchangeably.  

It is not hard to see that upon the above identification, we can define operations $e_r$ and $f_r$ on a diagram $D$ such that $(e_r(D))'=e_r^-(D')$, and similarly for $f_r$. Indeed, given a diagram $D$, we scan its cells in rows $r$ and $r+1$ columnwise from left to right and from top to bottom. We mark every cell in row $r$ (resp. $r+1$) as a left (resp. right) paranthesis, and match parantheses as usual. We define $e_r(D)$ to be the diagram obtained by moving the rightmost unpaired cell in row $r+1$ up by one position, if such a cell exists; otherwise, we define it to be the empty diagram $\emptyset$. Similarly, we define $f_r(D)$ to be the diagram obtained by moving the leftmost unpaired cell in row $r$ down by one position, if such a cell exists; otherwise, we define it to be the empty diagram. Furthermore, we can also define $\sigma_r(D)$ in the obvious way. Note that we can now restate the identity in Proposition \ref{plusminus} simply as 
\begin{equation}\label{erd}e_r^+(D)=e_r(D)\,,\end{equation}
and similarly for $f_r$.

\begin{exa}{\rm Here is an example of $r$-pairing on two successive rows of a diagram. The cells indicated with the same symbol except those marked with a circle are $r$-paired. The moves corresponding to the operations $e_r$ and $f_r$ are indicated by arrows.

\[
\begin{array}{c}
\mbox{\psfig{file=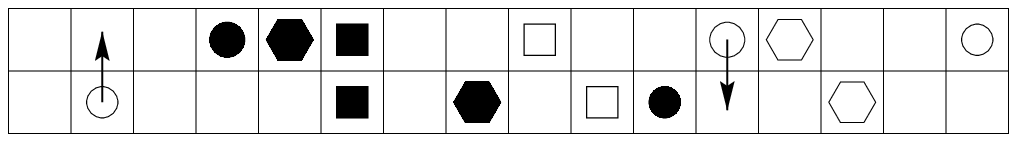}}
\end{array}
\]
}\end{exa}


We now list the properties of the operations $e_r$ and $f_r$ on diagrams. Also note that they define a directed graph $\widehat{\Gamma}$ on the set of all diagrams. 

\begin{thm}\label{erfrsr2} Let $h$ denote anyone of the operations $e_r$, $f_r$, $\sigma_r$, and let $D$ be a diagram.\\
\indent {\rm (1)} Assuming that $h(D)\ne \emptyset$, we have $P(h(D))=P(D)$.\\
\indent {\rm (2)} We have that $Q(h(D))=h(Q(D))$. \\
\indent {\rm (3)} The connected components of $\widehat{\Gamma}$ consist of all the diagrams with a fixed $P$-tableau. Two such components are isomorphic if and only if they are indexed by SSYT of the same shape.
\end{thm}

\begin{proof} The first property follows immediately from (\ref{erd}). The second property follows by combining (\ref{symm}) and Theorem \ref{erfrsr1} (2). The third property follows from the previous ones and Theorem \ref{erfrsr1} (4). 
\end{proof}

We now turn to the case of the nilplactic monoid, defined by the {\em Coxeter-Knuth} (or {\em nilplactic}) {\em equivalence} on reduced words for a permutation; this is the transitive closure of the relations
\begin{align*}
\ldots ikj\ldots\,&\sim_{CK}\,\ldots kij\ldots\\
\ldots jik\ldots\,&\sim_{CK}\,\ldots jki\ldots\\
\ldots i(i+1)i\ldots\,&\sim_{CK}\,\ldots (i+1)i(i+1)\ldots
\end{align*}
where $i<j<k$. We again consider biwords $w=\binom{u}{v}$ with $w^{\rm rev}$ in antilexicographic order, but now, throughout the remaining part of this section, $v$ is a reduced word (for some permutation). We can again associate with $w$ a pair of SSYT of conjugate shapes $(\widetilde{P}(w),\widetilde{Q}(w))$, with the only difference that $\widetilde{P}(w)$ is obtained by applying the {\em Edelman-Greene insertion} procedure \cite{eagbt} to $v$. Recall that, in a similar way to Knuth equivalence, two reduced words are Coxeter-Knuth equivalent if and only if they have the same Edelman-Greene insertion tableau. As mentioned in \cite{lasdcg}, we can define the operations $\widetilde{e}_r^+$ and $\widetilde{f}_r^+$ by using nilplactic jeu de taquin on two columns (denoted by $\widetilde{{\rm jdt}}$) and its inverse (denoted by $\widetilde{{\rm jdt}}^{-1}$). However, due to the lack of symmetry for the current type of biwords, the only way to define an analog of the operations $e_r^-$ and $f_r^-$ is via an analog of the diagram representation of a biword. The corresponding construction below is new.

We will identify the biword $w=\binom{u}{v}=\binom{u_1\ldots u_k}{v_1\ldots v_k}$ considered above with a line diagram $R$ similar to an rc-graph, which we call {\em stable rc-graph} (for reasons having to do with the theory of {\em stable Schubert polynomials}, also known as {\em Stanley symmetric functions}, cf. e.g. \cite{fasspn,manfsp}). We construct this as follows. Let $m$ be the last (largest) entry of $u$, and let $n-1$ be the largest entry of $v$ (thus, $v$ is a reduced word for some permutation $\pi_1\ldots\pi_n$ in $S_n$). The stable rc-graph $R$ will be the subset of $[m]\times[n+m]$ consisting of the pairs $(u_i,v_i+m-u_i+1)$, called {\em crosses}, which correspond to the biletters $\binom{u_i}{v_i}$ of $w$ (we let $[n]:=\{1,\ldots,n\}$ throughout this paper). This construction is best understood by considering the corresponding graphical representation, which is similar to that of an rc-graph. We draw $n$ lines going up and to the right, such that the $i$th line starts at position $(m,i)$ and ends at position $(1,\pi_i+m)$; as for rc-graphs, we use the matrix numbering of the positions. The rule for constructing the line diagram is the same as for rc-graphs (see the previous section). Given $R$ identified with $w=\binom{u}{v}$, we use the notation ${\rm scomp}(R)$ for $u$, and ${\rm red}(R)$ for $v$. We call the former a {\em semicompatible sequence}, since it satisfies only the first compatibility condition (\ref{comp1}). We also denote by ${\rm red}_i(R)$ the subword of ${\rm red}(R)$ corresponding to the $i$th row of $R$, that is, to entries $i$ in ${\rm scomp}(R)$. Note that when $u$ is a compatible sequence to $v$, namely $u_i\le v_i$ for all $i$, cf. (\ref{comp2}), then $R$ is essentially an rc-graph. 

\begin{exa}{\rm Here is an example of a stable rc-graph. It corresponds to the biword
\[\left(\begin{array}{cccccccccc}1&1&3&3&4&4&4&5&5&5\\3&1&4&2&3&2&1&5&4&2\end{array}\right)\,.\]
Thus $m=5$, $n=6$, and the corresponding permutation is $542613$.
\[
\begin{array}{c}
\mbox{\psfig{file=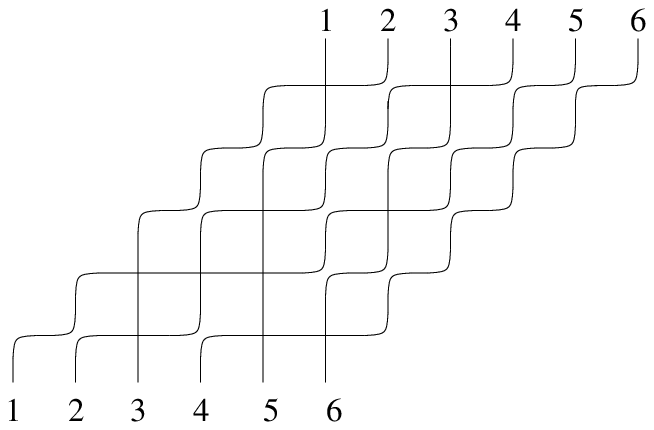}}
\end{array}
\]
}\end{exa}

The analogs $\widetilde{e}_r$, $\widetilde{f}_r$ of the operations $e_r^-$, $f_r^-$ above are now defined for stable rc-graphs $R$. We consider the crosses in rows $r$ and $r+1$ of $R$, and perform an $r$-pairing procedure identical to the one for diagrams. In order to define $\widetilde{e}_r(R)$, we perform an inverse chute move to the rightmost unpaired cross in row $r+1$, if such a cross exists. Similarly, we define $\widetilde{f}_r(R)$ by performing a chute move to the leftmost unpaired cross in row $r$, if such a cross exists. Furthermore, one can also define $\widetilde{\sigma}_r(R)$ by using a sequence of chute or inverse chute moves. Note that different moves on rc-graphs of the same nature were defined and studied recently in \cite{milmrc}.

\begin{rem}\label{compact}{\rm 
{\rm (1)} If crosses $(r,a)$ and $(r+1,b)$ form an $r$-pair, then all crosses $(r,c)$ and $(r+1,c)$ with $a<c<b$ are $r$-paired (with some other crosses in such positions). Hence rows $r$ and $r+1$ in an rc-graph consist of a sequence of blocks of $r$-paired crosses shuffled with a sequence of crosses in row $r+1$ followed by a sequence of crosses in row $r$.\\
\indent {\rm (2)} The inverse chute move in the definition of $\widetilde{e}_r$ is always possible if the corresponding unpaired cross exists. Indeed, the only way in which there is no inverse chute move applied to a cross $(r+1,a)$ is the one indicated in the figure below; namely, there is $b$ such that the rc-graph contains the crosses $(r,c)$ for $a<c<b$, and $(r+1,c)$ for $a\le c\le b$, but does not contain the cross $(r,b)$.
\[
\begin{array}{c}
\mbox{\psfig{file=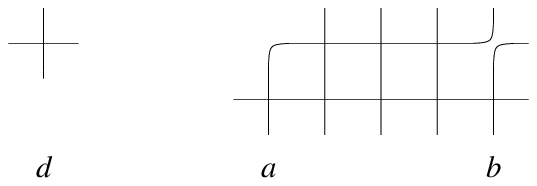}}
\end{array}
\]
Now the cross $(r+1,b)$ must be $r$-paired with a cross $(r,d)$, where $d<a$ (indeed, the cross $(r+1,a)$ is the righmost unpaired cross in row $r+1$). But this is a contradiction to the fact that the cross $(r+1,a)$ is unpaired, by the previous remark. There is a similar argument involving the operation $\widetilde{f}_r$.\\
\indent {\rm (3)} Assuming that $\widetilde{e}_r(R)\ne\emptyset$, it is not hard to see, based on the first remark, that $R$ and $\widetilde{e}_r(R)$ have the same $r$-pairing. Similarly for $\widetilde{f}_r$ and $\widetilde{\sigma}_r$.
}
\end{rem}

\begin{exa}\label{expair}{\rm Assume that rows $r$ and $r+1$ in a stable rc-graph are as in the figure below.
\[
\begin{array}{c}
\mbox{\psfig{file=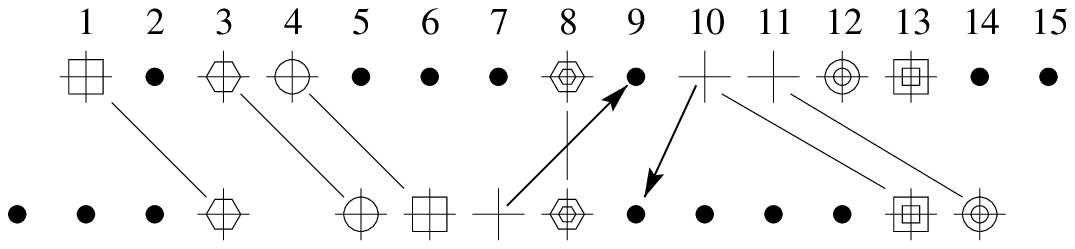}}
\end{array}
\]

The reduced word corresponding to the stable rc-graph contains the subword 
\[(13,12,11,10,8,4,3,1,15,14,9,8,7,6,4)\,.\]
The crosses marked with the same symbol are $r$-paired. The matching of certain crosses indicated by continuous lines is referred to in the proof of Proposition \ref{samep}. The moves corresponding to the operations $\widetilde{e}_r$ and $\widetilde{f}_r$ are indicated by arrows. 
}
\end{exa}

The properties of the operations $\widetilde{e}_r$, $\widetilde{f}_r$ on stable rc-graphs/biwords are easily deduced from those of $e_r$, $f_r$ on diagrams/biwords using the plactification map in \cite{raspla}. This is a map $\phi$ from reduced words to words which takes Coxeter-Knuth equivalence to Knuth equivalence. It is defined recursively by
\begin{equation}\label{defpl}\phi(\emptyset)=\emptyset\,,\;\;\;\;\;\phi(ru)=r\sigma_r(\phi(u))\,,\end{equation}
where $ru$ is a reduced word whose first entry is $r$. We summarize the properties of $\phi$ below.

\begin{thm}\label{pla}\cite{raspla} {\rm (1)} A reduced word $u$ is the row (resp. column) word of some skew tableau if and only if $\phi(u)$ is the row (resp. column) word of a skew tableau of the same shape.\\
\indent {\rm (2)} If $u$ and $\widetilde{u}$ are Coxeter-Knuth equivalent words, then $\phi(u)$ and $\phi(\widetilde{u})$ are Knuth equivalent. In other words, $\phi(\widetilde{P}(u))=P(\phi(u))$. In fact, more is true; namely, given reduced words $u'u u''$ and $u'\widetilde{u} u''$ with $u$ and $\widetilde{u}$ Coxeter-Knuth equivalent, we have $\phi(u'u u'')=v'v v''$ and $\phi(u'\widetilde{u} u'')=v'\widetilde{v} v''$, where $v$ and $\widetilde{v}$ are Knuth equivalent. In addition, the plactification map preserves the lengths of the corresponding subwords.\\
\indent {\rm (3)} We have $Q(\phi(u))=\widetilde{Q}(u)$ for all reduced words $u$.
\end{thm}

Note that the second part of property (2) is not explicitly stated in \cite{raspla}, but easily follows using the argument in the proof of Lemma 20 in that paper. In fact, this argument applies to the special case when $u$ and $\widetilde{u}$ are three letter subwords related by a Coxeter-Knuth transformation; the general case then easily follows from here. We can now state the following analog of Lemma \ref{uniqueright}; it is a straightforward consequence of this Lemma and of Theorem \ref{pla} (1)-(2).

\begin{lem}\label{uniqueright0}
Given two columns (decreasing words) $u,\,v$ of lengths $s$ and $t$, there is at most one pair of columns of lengths $s+1$ and $t-1$ in the nilplactic class of the word $uv$. If such a pair exists, then it is precisely $\widetilde{{\rm jdt}}(u,v)$. Similarly for $\widetilde{{\rm jdt}}^{-1}(u,v)$. 
\end{lem}

We now show that the two types of operations (giving left and right strings) on stable rc-graphs coincide. To this end, we first show that the operations $\widetilde{e}_r$ and $\widetilde{f}_r$ preserve the $\widetilde{P}$-tableau.

\begin{prop}\label{samep} Let $R$ be a stable rc-graph. If $\widetilde{e}_r(R)\ne\emptyset$, then $\widetilde{P}(\widetilde{e}_r(R))=\widetilde{P}(R)$. Otherwise, $\widetilde{{\rm jdt}}({\rm red}_r(R),{\rm red}_{r+1}(R))=\emptyset$. Similar statements hold for $\widetilde{f}_r$ and $\widetilde{\sigma}_r$.
\end{prop}
\begin{proof}
We can assume that $R$ has crosses only in rows $r$ and $r+1$. We examine the Edelman-Greene insertion procedure applied to the word $({\rm red}(R))^{\rm rev}$. It is well-known (cf. \cite{eagbt} Corollary 7.22) that the result of the insertion is $\widetilde{P}^T(R)$. Since the mentioned word is the concatenation of two rows (increasing words) $u,u'$, the insertion procedure can be described as follows. Scan the crosses in row $r$ from left to right, and match each one with the first (left-to-right) nonmatched cross in row $r+1$ in the same or a greater column. The matched crosses correspond to the pairs consisting of a letter in $u'$ to be inserted and the letter in $u$ it bumps; this still holds when the special bumping rule in the Edelman-Greene insertion \cite{eagbt} applies. Hence, the second row of $\widetilde{P}^T(R)$ consists of the letters of $u$ corresponding to the matched crosses in row $r+1$. The first row can now be easily found, based only on the permutation corresponding to $R$. Now let us perform the $r$-pairing procedure described above. Based on Remark \ref{compact} (1), we deduce that the $r$-paired crosses in row $r+1$ are precisely the previously matched crosses in this row (see the figure referring to Example \ref{expair}, where the matching is indicated by continuous lines). Hence $\widetilde{P}(R)$ is determined by the $r$-paired crosses in row $r+1$ and the permutation corresponding to $R$. The first part of the Proposition now follows by Remark \ref{compact} (3). The second part also follows; indeed, according to the remarks above, if all crosses in row $r+1$ are $r$-paired, then ${\rm red}(R)$ is precisely the column word of $\widetilde{P}(R)$. A similar argument holds for $\widetilde{f}_r$ and $\widetilde{\sigma}_r$.
\end{proof}

Combining Proposition \ref{samep} with Lemma \ref{uniqueright0}, we immediately deduce the analog of (\ref{erd}), namely Proposition \ref{plusminus0} below. Note that we could have used the same type of proof in the plactic case. However, we preferred the more conceptual proof based on the symmetry property of the RSK-correspondence, which does not work in the nilplactic case. 

\begin{prop}\label{plusminus0} For any stable rc-graph $R$, we have $\widetilde{e}_r^+(R)=\widetilde{e}_r(R)$. 
\end{prop}

We now define $\phi$ on a stable rc-graph $R$ to be the diagram identified with the biword
\begin{equation}\label{defphi}\phi(R):=\binom{{\rm scomp}(R)}{\phi({\rm red}(R))}\,.\end{equation}
By Theorem \ref{pla} (1), we can see that the biword $(\phi(R))^{\rm rev}$ is in antilexicographic order. In fact, more can be deduced using the same argument. To this end, let us first define the shape of a biword $w=\binom{u}{v}$ with $w^{\rm rev}$ in antilexicographic order, which is denoted by ${\rm shape}(w)$; this is the shape of the skew tableau obtained by placing the columns ${\rm col}_i(w)$ side by side, with a maximum overlap between every two consecutive columns, as indicated at the beginning of the section. We can see that ${\rm shape}(\phi(R))={\rm shape}(R)$ for any stable rc-graph $R$.

The following result will be used several times below.

\begin{prop}\label{commphi}
For any stable rc-graph $R$, we have $\phi(\widetilde{e}_r(R))={e}_r(\phi(R))$; similarly for $\widetilde{f}_r$ and $\widetilde{\sigma}_r$. 
\end{prop}
\begin{proof} This follows from Proposition \ref{plusminus0}, Lemma \ref{uniqueright0}, and Theorem \ref{pla} (1)-(2). Note that the strong form of Theorem \ref{pla} (2) is needed. Also note that $\widetilde{e}_r^+(R)=\emptyset$ if and only if $\widetilde{e}_r^+(\phi(R))=\emptyset$ because of the fact that ${\rm shape}(\phi(R))={\rm shape}(R)$. A similar argument holds for $\widetilde{f}_r$ and $\widetilde{\sigma}_r$.
\end{proof}

Finally, we can prove the following analog of Theorem \ref{erfrsr2} (2).

\begin{prop}\label{erfr0} For any stable rc-graph $R$, we have $\widetilde{Q}(\widetilde{e}_r(R))=e_r(\widetilde{Q}(R))$; similarly for $\widetilde{f}_r$ and $\widetilde{\sigma}_r$.
\end{prop}
\begin{proof} This is a consequence of the commutativity of the diagram below; here ${\mathcal S\mathcal R}$ denotes the set of all stable rc-graphs, ${\mathcal Diag}$ the set of all diagrams, and ${\mathcal T\!ab}$ the set of SSYT.

\[
\begin{diagram}
\node[2]{{\mathcal S\mathcal R}}\arrow[2]{e,t}{\widetilde{e}_r,\widetilde{f}_r,\widetilde{\sigma}_r}
\arrow{s,r}{\phi}\arrow[1]{ssw,l}{\widetilde{Q}}
\node[2]{{\mathcal S\mathcal R}}\arrow{s,l}{\phi}\arrow[1]{sse,t}{\widetilde{Q}}\\
\node[2]{{\mathcal Diag}}\arrow[2]{e,t}{e_r,f_r,\sigma_r}\arrow{sw,b}{Q}
\node[2]{{\mathcal Diag}}\arrow{se,b}{Q}\\
\node{{\mathcal T\!ab}}\arrow[4]{e,t}{e_r,f_r,\sigma_r}\node[4]{{\mathcal T\!ab}}
\end{diagram}
\]

The two triangles on the sides commute by Theorem \ref{pla} (3). The rectangle above commutes by Proposition \ref{commphi}. The trapezoid below commutes by Theorem \ref{erfrsr2} (2).
\end{proof}

\section{Key Polynomials}

In this section we investigate in more detail the expression of a Schubert polynomial in terms of key polynomials given in \cite{lastns}. In particular, we relate the combinatorics of key polynomials in terms of SSYT to the combinatorics of Schubert polynomials in terms of rc-graphs. 

Key polynomials (for the symmetric group) are certain polynomials in $\bZ[x_1,x_2,\ldots]$ indexed by compositions $\alpha$, which interpolate between monomials and Schur polynomials. From a representation theoretic point of view, they are certain ``partial'' characters of the irreducible representations of $GL_n(\bC)$. More precisely, the key polynomial $\kappa_\alpha(X)$ is defined by 
\[\kappa_\alpha(X):=\pi_{u(\alpha)}(x^{\lambda(\alpha)})\,.\]
Here $x^\beta:=x_1^{\beta_1}x_2^{\beta_2}\ldots$ for any composition $\beta=(\beta_1,\beta_2,\ldots)$, and $u(\alpha)$ is the permutation of minimum length which realizes the rearrangement $\lambda(\alpha)$ of $\alpha$ in weakly decreasing order via the right action $\alpha\,w:=(\alpha_{w_1},\alpha_{w_2},\ldots)$. On the other hand, $\pi_w$ is the isobaric divided difference operator on $\bZ[x_1,x_2,\ldots]$, defined by
\[\pi_w:=\pi_{a_1}\ldots\pi_{a_{l(w)}}\,,\spa\spa\spa\spa\pi_r:=\partial_r x_r\,;\]
here $\partial_r$ are the divided difference operators introduced in Section 1, $a_1\ldots a_{l(w)}$ is any reduced word for $w$, $x_r$ denotes multiplication by this variable, and a permutation $w$ acts on $\bZ[x_1,x_2,\ldots]$ by $w\,x_r:=x_{w_r}$. 

In order to study the combinatorics of key polynomials, we need the concepts of left and right keys of a tableau.  The {\em left key} $K_{-}(T)$ (resp. {\em right key} $K^+(T)$) of a SSYT $T$ (cf. \cite{lasksb}) is the tableau of the same shape as $T$ whose $j$th column is the first (resp. last) column of any skew tableau in the plactic class of $T$ with the following properties: (1) its sequence of column lengths is a rearrangement of the corresponding sequence for $T$; (2) its first (resp. last) column has the same length as the $j$th column of $T$. It is not obvious that this is a correct definition. An easy way to find the $j$th column of $K_{-}(T)$ is to perform the reverse column insertion procedure to the entries in the $j$th column of $T$, successively. Alternatively, we can use jeu de taquin on pairs of successive columns to interchange the lengths of columns in $T$. Similar procedures apply to $K_+(T)$. On the other hand, by working in the nilplactic monoid instead of the plactic one, we can define the {\em left nil key} $\widetilde{K}_{-}(T)$ (resp. {\em right nil key} $\widetilde{K}^+(T)$) of a SSYT $T$ whose row (or column) word is a reduced word. All four types of tableaux defined above have the property that the entries in any column except the first form a subset of the entries in the previous column. Such a tableau is called a {\em key}. There is an obvious bijection between compositions and keys given by $\alpha\mapsto{\rm key}(\alpha)$, where ${\rm key}(\alpha)$ is the SSYT of shape $\lambda(\alpha)$ whose first $\alpha_j$ columns contain the letter $j$, for all $j$. The inverse map is given by $T\mapsto{\rm content}(T)$, where ${\rm content}(T)$ is the composition whose $j$th part is the number of entries equal to $j$ in $T$. Given a composition $\alpha$, we will use the following notation
\[\mathcal T\!ab(\alpha):=\{\,T\mbox{ SSYT of shape }\lambda(\alpha)\,,\;K_+(T)\le{\rm key}(\alpha)\,\}\,;\]
here the inequality means entrywise comparison. We can now state the following combinatorial construction of key polynomials in \cite{lasksb}.

\begin{thm}\label{keypol}\cite{lasksb} For any composition $\alpha$, we have
\[\kappa_\alpha(X)=\sum_{T\in\mathcal T\!ab(\alpha)}x^T\,,\]
where $x^T:=x^{{\rm content}(T)}$.
\end{thm}

In order to prove this result, Lascoux introduced combinatorial versions of the operators $\pi_r$, which are set-valued maps from SSYT to SSYT. More precisely, given a tableau $T$, one defines $\pi_r(T)$ to be the set of tableaux $\{T,f_r(T),\ldots,f_r^s(T)\}$ if the subword of unpaired $r$'s and $r+1$'s in the column word of $T$ is of the form $r^s$; otherwise it is defined to be the empty set. More generally, given a set-valued map $f$ from a set $X$ to $X$ and a multiset $M$ with elements from $X$, one defines $f(M)$ to be the {\em multiset} union of $f(a)$ for each $a$ in $M$, taken with its corresponding multiplicity; in particular $f(\emptyset)=\emptyset$. Given a finite sequence of such set-valued maps, one can define the value of their composition on a multiset in the obvious way, using the above definition. Note that, in general, the composition $\pi_{a_1}\ldots\pi_{a_p}(T)$ gives different results for different reduced words $a_1\ldots a_p$ for a given permutation. We can now state the following result, which is contained in the proof of Theorem \ref{keypol} given in \cite{lasksb}. Recall the notation $U(\lambda)$ for the tableau of shape $\lambda$ with all entries in row $i$ equal to $i$.

\begin{thm}\cite{lasksb}\label{combkey}
Let $\alpha$ be a composition, and let $a_1\ldots a_p$ be {\em any} reduced word for $u(\alpha)$. Then we have
\[\mathcal T\!ab(\alpha)=\pi_{a_1}\ldots\pi_{a_p}(U(\lambda(\alpha)))\,;\]
here composition of set-valued maps is understood in the sense mentioned above. 
\end{thm}

\begin{rema}\label{referee}{\rm  Littelmann has a vast generalization of Theorem \ref{combkey} \cite{litlrr,litpro}. His paths can be generated for any {\em Demazure module} (over a symmetrizable Kac-Moody algebra) in the same kind of fashion. The idea of left and right key is a type $A$ translation of the notion of a {\em canonical lift} in the standard monomial theory of Lakshmibai, Musili, and Seshadri \cite{lassmt1}. The {\em initial direction} of a Lakshmibai-Seshadri path (which is a special case of a Littelmann path) is the extremal weight which gives the ``right key'' of the path; the {\em final direction} is the ``left key''. }
\end{rema}

Consider now the set ${\cal R}(w)$ of rc-graphs corresponding to a permutation $w$. We identify an rc-graph $R$ with the pair $\binom{{\rm comp}(R)}{{\rm red}(R)}$. We associate to $R$ the pair of SSYT of conjugate shapes $(\widetilde{P}(R),\widetilde{Q}(R))$ (which determine $R$), as in Section 3; note that these tableaux were already considered in \cite{bjsscp}. We denote by ${\mathcal P}(w)$ the set of row and column strict tableaux indexing the nilplactic classes of reduced words for $w$; more explicitly, the elements of ${\mathcal P}(w)$ correspond to those reduced words for $w$ which are column (or row) words for a tableau. We also use the following notation:
\[{\mathcal R}(P):=\setp{R\in{\mathcal R}(w)}{\widetilde{P}(R)=P}\,,\;\;\;\;\;{\mathcal Q}(P):=\setp{\widetilde{Q}(R)}{R\in{\mathcal R}(P)}\,.\]
The sets ${\mathcal R}(P)$ were called {\em crystals} in \cite{lasdcg}. They partition the set of rc-graphs ${\mathcal R}(w)$ into blocks indexed by the tableaux in ${\mathcal P}(w)$. The map $R\mapsto \widetilde{Q}(R)$ is a bijection between ${\mathcal R}(P)$ and ${\mathcal Q}(P)$; furthermore, $x^R=x^{\widetilde{Q}(R)}$, so it makes sense to call the bijection monomial preserving. 

The question was raised in \cite{bjsscp} to characterize the tableaux in ${\mathcal Q}(P)$. It can be proved using the results and techniques in \cite{raskpf}, cf. also \cite{raspla}, that a stable rc-graph $R$ is an rc-graph (that is, the second compatibility condition (\ref{comp2}) is satisfied) if and only if
\begin{equation}\label{keycond}
K_+(\widetilde{Q}(R))\le\widetilde{K}_-(\widetilde{P}^T(R))\,.
\end{equation}
Essentially, the proof is a careful examination of the way in which the mentioned keys are constructed via reverse column insertion. This procedure allows one to translate the compatibility condition (\ref{comp2}) into the language of tableaux. It follows that
\begin{equation}\label{charq}
{\mathcal Q}(P)=\setp{Q}{K_+(Q)\le\widetilde{K}_-(P^T)}=\mathcal T\!ab({\rm content}(\widetilde{K}_-(P^T)))\,.
\end{equation}
Hence, the set ${\mathcal Q}(P)$, which is in a monomial preserving bijection with ${\mathcal R}(P)$, underlies the combinatorial construction of a certain key polynomial (by Theorem \ref{keypol}). The following results about key polynomials now immediately follow: (1) the expansion of a Schubert polynomial as a positive sum of key polynomials; (2) a combinatorial formula for key polynomials in terms of rc-graphs. 

\begin{thm}\label{las1}\cite{lastns,raskpf} We have that
\[\fs_w(X)=\sum_{P\in{\mathcal P}(w)}\kappa_{{\rm content}(\widetilde{K}_{-}(P^T))}(X)\,.\]
Thus, we have the following combinatorial formula for Schubert polynomials:
\[\fs_w(X)=\sum_{P\in{\mathcal P}(w)}\;\,\sum_{Q\in\mathcal T\!ab({\rm content}(\widetilde{K}_{-}(P^T)))}x^Q\,.\]
\end{thm}

\begin{thm}\label{las2}\cite{raskpf} Given a composition $\alpha$, and any SSYT $P$ in ${\mathcal P}(w)$ for some permutation $w$, with ${\rm content}(\widetilde{K}_{-}(P^T))=\alpha$, we have
\[\kappa_\alpha(X)=\sum_{R\in{\mathcal R}(P)} x^R\,.\]
\end{thm}

We conclude this section with an investigation of the structure of the crystals ${\mathcal R}(P)$. In particular, we are interested in an efficient procedure to generate a crystal. 

\begin{rema}\label{oponrc}{\rm
According to Proposition \ref{samep}, for every $R$ in ${\mathcal R}(P)$, we have that $\widetilde{e}_r(R)$ is in ${\mathcal R}(P)$ too, unless it is $\emptyset$. Furthermore, according to Proposition \ref{erfr0}, we have that ${e}_r(\widetilde{Q}(R))$ is in ${\mathcal Q}(P)$, unless it is $\emptyset$, and, in fact, it coincides with $\widetilde{Q}(\widetilde{e}_r(R))$. Hence, there is a perfect correspondence between the action of $\widetilde{e}_r$ on ${\mathcal R}(P)$ and that of the usual coplactic operations on the tableaux in ${\mathcal Q}(P)$. On the other hand, it is possible that $\widetilde{f}_r(R)\ne\emptyset$ is a stable rc-graph which is not an rc-graph; indeed, assume, for instance, that $R$ has only one cross in rows $r$, in position $(r,1)$, and no crosses in row $r+1$.}
\end{rema} 

The set ${\mathcal Q}(P)$ has a natural partial order by reverse entrywise comparison: $Q\preceq Q'$ if and only if $Q\ge Q'$. This can be transferred to ${\mathcal R}(P)$ via the bijection $R\mapsto\widetilde{Q}(R)$. We have the following result.

\begin{prop}\label{thmtwo}
{\rm (1)} The poset $({\mathcal R}(P),\preceq)$ has a maximum $R_{\rm top}(P)$ and a minimum $R_{\rm bot}(P)$. These are determined by
\[\widetilde{Q}(R_{\rm top}(P))=U({\rm shape}(P^T))\,,\;\;\;\;\;\;\;\;\;\;\widetilde{Q}(R_{\rm bot}(P))=\widetilde{K}_{-}(P^T)\,.\]
\indent {\rm (2)} We have that $R\prec\widetilde{e}_r(R)$ whenever $\widetilde{e}_r(R)\ne\emptyset$. Furthermore, any rc-graph in ${\mathcal R}(P)$ can be obtained from $R_{\rm top}(P)$ by successively applying operations $\widetilde{f}_r$, but not all rc-graphs in this crystal can be obtained from $R_{\rm bot}(P)$ by applying operations $\widetilde{e}_r$. 
\end{prop} 
\begin{proof}
Most of the proof relies on Remark \ref{oponrc}. The first part of (2) is clear. Note that it is possible to successively apply coplactic operations $e_r$ to a SSYT of shape $\lambda$ until we reach $U(\lambda)$. Since $\widetilde{e}_r$ applied to an rc-graph always gives another rc-graph or $\emptyset$, while $\widetilde{f}_r$ might give a stable rc-graph which is not an rc-graph, the second part of (2) also follows (in fact, it is easy to construct an example showing that $R_{\rm bot}(P)$ might not be reachable from certain rc-graphs in the crystal via operations $\widetilde{f}_r$). The assertion involving $R_{\rm top}(P)$ in (1) follows from the previous arguments. The assertion involving $R_{\rm bot}(P)$ follows from the characterization (\ref{charq}) of ${\mathcal Q}(P)$ and the fact that $K_+(T)\ge T$ for any SSYT $T$ (this can be easily seen using one of the two constructions of the right key mentioned above). 
\end{proof}

Clearly, the rc-graphs $R_{\rm top}(w)$ and $R_{\rm bot}(w)$ are among the rc-graphs $R_{\rm top}(P)$ and $R_{\rm bot}(P)$. Notice that the second part of Proposition \ref{thmtwo} gives an algorithm for generating the crystal indexed by $P$, starting from $R_{\rm top}(P)$. The similar algorithm based on nilplactic jeu de taquin on two columns was illustrated in \cite{lasdcg}. We will now present a more efficient algorithm, in which each rc-graph in the crystal is obtained precisely once. To this end, we need to define combinatorial versions of the operators $\pi_r$ for rc-graphs, which will be denoted by $\widetilde{\pi}_r$ (to make the notation consistent). This can be done in a similar way to the definition of such operators on SSYT. Note that we set $\widetilde{\pi}_r(R):=\emptyset$ unless all the unpaired crosses in rows $r$ and $r+1$ of $R$ lie in the the former row. We can now transfer the result of Theorem \ref{combkey} from ${\mathcal Q}(P)$ to ${\mathcal R}(P)$ via the bijection $\widetilde{Q}(R)\mapsto R$. Indeed, due to Proposition \ref{erfr0}, the action of the operators $\widetilde{\pi}_r$ on rc-graphs corresponds to that of $\pi_r$ on tableaux via the mentioned bijection.

\begin{thm}\label{gencrystal}
Let $P$ be a tableau in ${\mathcal P}(w)$, and let $a_1\ldots a_p$ be {\em any} reduced word for $u({\rm content}(\widetilde{K}_-(P^T)))$. Then we have
\[{\mathcal R}(P)=\widetilde{\pi}_{a_1}\ldots\widetilde{\pi}_{a_p}(R_{\rm top}(P))\,;\]
here composition of set-valued maps is understood in the sense mentioned before Theorem \ref{combkey}.
\end{thm}

Theorem \ref{gencrystal} gives an efficient algorithm for generating the crystal ${\mathcal R}(P)$. This can be illustrated by constructing a tree whose root corresponds to $R_{\rm top}(P)$, and whose edges correspond to applying the operators $\widetilde{\pi}_r$ on rc-graphs. Note that $R_{\rm top}(P)$ can be easily constructed by reverse Edelman-Greene insertion, using the fact that $\widetilde{Q}(R_{\rm top}(P))=U({\rm shape}(P^T))$, cf. Proposition \ref{thmtwo} (1). To make things more explicit, we consider the example in \cite{raskpf}.

\begin{exa}\label{excrys}{\rm Let $w=21543$. There are three crystals corresponding to this permutation. We choose the one indexed by the tableau $P$ indicated below, together with $\widetilde{K}_-(P^T)$:
\[\begin{array}{lllll}
1&3&4\;\;\;\;\;\;\;&1&3\\
4&\,&\,&3\\
\,&\,&\,&4
\end{array}
\]
The composition $\alpha:={\rm content}(\widetilde{K}_-(P^T))$ is $(1,0,2,1)$, and the permutation $u(\alpha)$ is $3142$. We consider the reduced word $213$ for it. The construction of the crystal ${\mathcal R}(P)$ is indicated in the appendix to the paper.
}
\end{exa}

\begin{rem}\label{vexrem}
{\rm (1) There is a version of Theorem \ref{las2} in terms of words rather than reduced words (see \cite{raskpf} Theorem 5 (2) and Remark \ref{impla} (2)). This can be obtained by using the plactification map mentioned in the previous section (cf. Remark \ref{impla} (2)).\\
\indent (2) The permutations $w$ for which there is a unique crystal are known as {\em vexillary}. The corresponding Schubert polynomials, which thus coincide with a certain key polynomial $\kappa_\alpha(X)$ (in fact, $\alpha={\rm code}(w)$), are known to be {\em flagged Schur polynomials} (see \cite{macnsp,manfsp,raskpf}). The flag condition is given by a certain sequence called the {\em flag} of $w$ and denoted $\phi(w)$, which is defined as the weakly increasing rearrangement of the numbers
\[d_i:=\min\,\setp{j}{j>i\mbox{ and } w_j<w_i}-1\,.\]
The usual proof of this fact is based on algebraic manipulations with divided difference operators. We are able to prove this fact bijectively, by showing that the map $R\mapsto \widetilde{Q}(R)$ is a bijection from ${\mathcal R}(w)$ to the set of SSYT of shape $\lambda({\rm code}(w))$ satisfying the flag condition. More precisely, we can show that the flag condition mentioned above is equivalent to the general inequality condition (\ref{keycond}) by using the combinatorics of vexillary permutations; in particular, we need the characterization of their code, as well as the construction of ${\rm code}(w)$ which inputs $\lambda({\rm code}(w))$ and $\phi(w)$, cf. \cite{macnsp} p. 13--16. The proof is somewhat technical, and therefore we decided to omit it. However, we note that this approach emphasizes the importance of condition (\ref{keycond}), which will be mentioned later in this paper, as well (in Remark \ref{impla} (2)).\\
\indent (3) More generally, given a Schubert polynomial which is a {\em flagged skew Schur polynomial} (see Theorem 26 in \cite{raskpf} for a sufficient condition), we can ask for a bijection, possibly monomial preserving, between the corresponding rc-graphs and skew tableaux satisfying the flag condition. See Remark \ref{balrem} (3) for the case of 321-avoiding permutations.\\
\indent (4) It is known that the crystal graph structure on tableaux is relevant to the representation theory of $GL_n(\bC)$. On the other hand, it was shown in \cite{raskpf} that key polynomials are characters of certain representations of the subgroup $B$ of $GL_n(\bC)$ consisting of upper triangular matrices. The corresponding modules were also constructed explicitly. The $r$-pairings and the operations $\widetilde{e}_r$, $\widetilde{f}_r$, $\widetilde{\pi}_r$ on rc-graphs considered above have a clear meaning in this context, but this is not presented here.
}
\end{rem}

\section{Splitting Schubert polynomials}

In this section, we use the operations $\widetilde{f}_r$ on rc-graphs to give a straightforward, purely combinatorial proof of a formula in \cite{bktspq}, which expresses Schubert polynomials as linear combinations of products of Schur polynomials in different sets of variables. This formula is one of the main results in \cite{bktspq}, and is derived in connection with {\em quiver varieties}. More precisely, it was shown that the coefficients in the mentioned formula are all {\em quiver coefficients}. Note that the derivation in \cite{bktspq} relies on a formula of Fomin and Greene, which gives a combinatorial interpretation for the coefficients in the expansion of stable Schubert polynomials (also known as Stanley symmetric functions) in the basis of Schur functions \cite{fagnsf}. The result of Fomin and Greene, which is obtained as an application of their theory of {\em noncommutative Schur functions}, is not needed in the proof below. 

Let us now recall the setup in \cite{bktspq}. First recall that a permutation $w$ is said to have a {\em descent} in position $i$ whenever $w_i>w_{i+1}$. We say that a set of integers $A$ is {\em compatible with} $w$ if all descent positions of $w$ are contained in $A$. The Schur polynomial indexed by a partition $\lambda$ in a set of variables $X$ is denoted by $s_\lambda(X)$. The column word of a SSYT $T$ is denoted by ${\rm col}(T)$. We let $U(\lambda,a)$ denote the SSYT of shape $\lambda$ with all entries in row $i$ equal to $a+i$. We say that a sequence of SSYT $(T_1,\ldots,T_k)$ is {\em strictly bounded below} by a sequence $(a_1,\ldots,a_k)$ if the entries of $T_i$ are strictly greater than $a_i$, for each $i$. 

\begin{thm}\cite{bktspq} Suppose that $w$ is a permutation in $S_n$ compatible with the set of integers $\{0=a_0< a_1<\ldots<a_m\}$. Then we have
\[\fs_w(X)=\sum_\lambda c_\lambda s_{\lambda^1}(X_1)\ldots s_{\lambda^m}(X_m)\,,\]
where $X_i=\{x_{a_{i-1}+1},\ldots,x_{a_i}\}$, and the sum is over all sequences of partitions $\lambda=(\lambda^1,\ldots,\lambda^m)$. The coefficients $c_\lambda$ are nonnegative integers; furthermore, they count the sequences of SSYT $(T_1,\ldots,T_m)$ strictly bounded below by $(0,a_1,\ldots,a_{m-1})$, such that the shape of $T_i$ is conjugate to $\lambda^i$ and ${\rm col}(T_1)\ldots{\rm col}(T_m)$ is a reduced word for $w$.
\end{thm}
\begin{proof}
Assume $a_m<n$. Given an rc-graph $R$ in ${\mathcal R}(w)$, we can clearly split it into a sequence of rc-graphs $(R_1,\ldots,R_m)$, where $R_k:=\setp{(i,j)\in R}{a_{k-1}+1\le i\le a_k}$, for each $k$. Note that $w(R_k)$ is not necessarily a permutation with a unique descent (that is, a {\em Grassmannian permutation}). We associate to $R_k$ the pair of SSYT $(\widetilde{P}(R_k),\widetilde{Q}(R_k))$, as in Section 3. It is easy to check that the above procedure defines a map from ${\mathcal R}(w)$ to sequences of pairs of SSYT of conjugate shapes $((P_1,Q_1),\ldots,(P_m,Q_m))$, which satisfy the following conditions:
\begin{enumerate}
\item $(P_1,\ldots,P_m)$ is strictly bounded below by $(0,a_1,\ldots,a_{m-1})$;
\item ${\rm col}(P_1)\ldots{\rm col}(P_m)$ is a reduced word for $w$;
\item the entries of $Q_k$, for $k=1,\ldots,m$, are between $a_{k-1}+1$ and $a_k$. 
\end{enumerate}
Clearly, this map is injective, so we only need to prove surjectivity.  

Consider a sequence of pairs of tableaux $((P_1,Q_1),\ldots,(P_m,Q_m))$ with the above properties. This determines a sequence of stable rc-graphs $(R_1,\ldots,R_m)$, which, in turn, determines a stable rc-graph for the permutation $w$, by condition (2) (just by concatenating the corresponding reduced words and compatible sequences, or, graphically, by stacking the relevant portions of $R_i$ on top of each other). It suffices to show that each $R_k$ is, in fact, an rc-graph. Let $\lambda^k$ be the shape of $Q_k$. It is not hard to see, based on the first two conditions above, that $(P_k,U(\lambda^k,a_{k-1}))$ determines an rc-graph $R_k^0$. We can get from the tableau $U(\lambda^k,a_{k-1})$ to $Q_k$ by succesively applying operations $f_r$, with $a_{k-1}+1\le r\le a_k-1$ (due to condition (3)). We will show that the same sequence of operations $\widetilde{f}_r$ applied to $R_k^0$ produces an rc-graph. Then we conclude that this rc-graph has to be $R_k$, based on Propositions \ref{samep} and \ref{erfr0} (which state the compatibility of the actions of $\widetilde{f}_r$ on rc-graphs  and of $f_r$ on the associated tableaux).

Now consider the above sequence of operations $\widetilde{f}_r$ applied to the concatenation of the rc-graphs $R_1^0,\ldots,R_m^0$, which is an rc-graph for $w$. Based on the properties of the operations $\widetilde{f}_r$ discussed above, it suffices to show that, when applying $\widetilde{f}_r$ to an rc-graph for $w$, we never obtain a stable rc-graph which is not an rc-graph. This could only happen if, at the beginning of rows $r$ and $r+1$, we have the configuration shown below, where the circled cross is the one moved by $\widetilde{f}_r$. 
\[
\begin{array}{c}
\mbox{\psfig{file=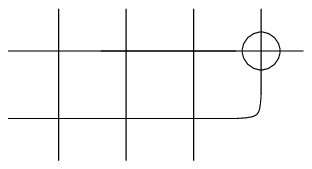}}
\end{array}
\]
But this means that the permutation $w$ has a descent in in position $r$, which is a contradiction with the compatibility of the sequence $\{0=a_0< a_1<\ldots<a_m\}$.
\end{proof}

Let us note, as it is done in \cite{bktspq}, that the classical construction of Schubert polynomials (\ref{schub1}) is equivalent to the special case of the above theorem corresponding to the set of integers $[n-1]$.

\section{Kohnert Diagrams}

Kohnert \cite{kohwpt} conjectured that Schubert polynomials can be constructed in terms of certain diagrams, currently known as {\em Kohnert diagrams}. These are subsets of $[n]\times[n]$ which can be obtained from the {\em diagram of a permutation} by certain moves (recall that $[n]:=\{1,\ldots,n\}$). The diagram $D(w)$ of the permutation $w$ is the set
\[D(w):=\setp{(i,w_j)}{i<j,\: w_i>w_j}\subseteq[n]\times[n]\,.\]
 A {\em Kohnert move} (or, simply, {\em K-move}) on a Kohnert diagram picks the rightmost cell in any row $i$, and moves it up to the first nonoccupied cell directly above it, say in row $j$; we indicate such a move by $i\rightarrow j$. Let ${\mathcal K}(w)$ denote the closure of $\{D(w)\}$ under all possible K-moves. Kohnert's conjecture says that
\begin{equation}\label{koh}\fs_w(X)=\sum_{K\in{\mathcal K}(w)}x^K\,,\end{equation}
where $x^K:=\prod_{(i,j)\in K} x_i$. 

Bergeron \cite{berccs} defined more complex moves on diagrams, called {\em B-moves}, and constructed the set of diagrams ${\mathcal B}(w)$ as the closure of $\{D(w)\}$ under all possible B-moves. 
He stated that
\begin{equation}\label{ber}\fs_w(X)=\sum_{B\in{\mathcal B}(w)}x^B\,.\end{equation}
He used Schubert induction, starting with the longest permutation $w_0^n$ in $S_n$, and successively applying combinatorial versions of the operators $\partial_r$ for diagrams. The definition of the latter can be easily rephrased using the operations $f_r$ on diagrams in Section 3. More precisely, if $w\ne w_0^n$, Bergeron considered $r:=\min\,\setp{i}{w_i<w_{i+1}}$. Then, given a diagram $B$ in ${\mathcal B}(ws_r)$, he considered $B':=B\setminus\{(r,w_r)\}$. Finally, he defined $\partial_r(B)$ to be the set of diagrams $\{B',f_r(B'),\ldots,f_r^t(B')\}$ if all the $t$ unpaired cells in the $r$-pairing for $B'$ are in row $r$; otherwise, $\partial_r(B)$ was defined to be the empty set. Bergeron stated the following combinatorial version of the identity $\fs_w(X)=\partial_r(\fs_{ws_r}(X))$ defining Schubert polynomials, which implies (\ref{ber}):
\begin{equation}\label{combb}
{\mathcal B}(w)=\partial_r({\mathcal B}(ws_r))\,;\end{equation}
the right-hand side is understood in the multiset sense discussed before Theorem \ref{combkey}. 

Winkel \cite{windrg} announced a proof of Kohnert's conjecture (\ref{koh}) and of the fact that ${\mathcal K}(w)={\mathcal B}(w)$, which was conjectured in \cite{berccs}. He thus obtained (\ref{ber}) as a Corollary. This proof is still by induction, and is still based on Bergeron's combinatorial $\partial_r$, but is very involved. The reason for this is the fact that the action of $\partial_r$ on diagrams does not have a close relationship with K-moves, although it does have a closer one with B-moves. Later, Winkel announced a shorter proof of just (\ref{koh}) in \cite{winspk}, based on {\em Monk's formula} for Schubert polynomials. But this proof is still not transparent, using complex manipulations with Kohnert diagrams. In fact, Winkel himself recognized that his shorter proof is ``done by an exhaustive study of cases, which establishes correctness without leading to a deeper understanding''. 

We attempted to find a more transparent proof of Kohnert's conjecture by constructing a bijection between Kohnert diagrams and rc-graphs based on relating K-moves and L-moves. Although these moves look very similar in nature, their relationship is complicated, so we failed in this attempt. 

One of the reasons for the difficulties in all the mentioned attempts is that, as opposed to rc-graphs, there is no intrinsic description of the diagrams in ${\mathcal B}(w)$ or ${\mathcal K}(w)$, just recursive constructions. We now investigate a new set of diagrams $\varPhi(w)$, which is in a monomial preserving bijection with ${\mathcal R}(w)$ (thus providing a combinatorial construction of $\fs_w(X)$), and which has both a recursive and a nonrecursive construction. Furthermore, we show that $\varPhi(w)$ coincides with ${\mathcal B}(w)$ and ${\mathcal K}(w)$, if we assume the results of Bergeron and Winkel mentioned above. The nonrecursive construction of $\varPhi(w)$ is very easy to state using the one-to-one map $\phi$ on rc-graphs defined in (\ref{defphi}), namely $\varPhi(w):=\phi({\mathcal R}(w))$. Clearly, $\phi$ is a monomial preserving bijection from ${\mathcal R}(w)$ to $\varPhi(w)$, that is, $x^{\phi(R)}=x^R$. 

\begin{rem}\label{impla}{\rm (1) In \cite{raspla}, the image of the plactification map on reduced words for $w$ is characterized in terms of {\em $D(w^{-1})$-peelable words}. Thus, we have the following nonrecursive description of $\varPhi(w)$:
\begin{align}\label{phiw}\varPhi(w)=\{\:D\subseteq[n]\times[n]\::\:&{\rm col}(D)\mbox{ is }D(w^{-1})\mbox{-peelable and }\\&{\rm row}(D)\le\phi^{-1}({\rm col}(D))\:\}\,, \nonumber\end{align}
where $\le$ means entrywise comparison. One might be able to rephrase the above inequality in an even more explicit way; essentially, we are asking: what do compatible sequences for $u$ become under the map $u\mapsto\phi(u)$? \\
\indent (2) The situation is nicer if we replace the map $u\mapsto\phi(u)$ in the previous remark with $u\mapsto(\phi(u^{\rm rev}))^{\rm rev}$. Indeed, it turns out that the set of compatible sequences is preserved under the latter map (compatible sequences for arbitrary words are defined in the same way as for reduced words). Here is an easy way to see this: we translate the compatibility condition into the language of keys as stated in (\ref{keycond}) (there is a similar statement for words instead of reduced words), and then we invoke the plactification map. We need, of course, the properties of this map listed in Theorem \ref{pla}, as well as the following property appearing in \cite{raspla}: $K_-(\phi(P))=\widetilde{K}_-(P)$, for SSYT $P$ whose column word is reduced. The mentioned fact implies another combinatorial construction for $\fs_w(X)$ (in terms of $D(w)$-peelable words), which was stated in \cite{raspla} as Theorem 31 (1). }
\end{rem}

We now turn to the recursive construction of $\varPhi(w)$. The idea is to first give a recursive construction of ${\mathcal R}(w)$ based on analogs $\widetilde{\partial}_r$ of the operators $\partial_r$ on diagrams defined above. Then we show that the latter correspond to the former via the bijection $\phi$. Given a permutation $w\ne w_0^n$ in $S_n$, we let $r$ be defined as above, namely $r:=\min\,\setp{i}{w_i<w_{i+1}}$. Then, given an rc-graph $R$ in ${\mathcal R}(ws_r)$, we proceed in the same way as for diagrams by letting $R':=R\setminus\{(r,w_r)\}$ (the rc-graph is viewed as a subset of $\bZ_{>0}\times\bZ_{>0}$ here). Finally, we define $\widetilde{\partial}_r(R)$ to be the set of rc-graphs $\{R',\widetilde{f}_r(R'),\ldots,\widetilde{f}_r^t(R')\}$ if all the $t$ unpaired crosses in the $r$-pairing for $R'$ are in row $r$; otherwise, it is defined to be the empty set. The fact that $R'$ is indeed an rc-graph (i.e., its lines have no multiple intersections) can be seen by examining the figure below, which shows the form of {\em any} rc-graph in ${\mathcal R}(ws_r)$; the circled cross is the one to be removed. On the other hand, the successive applications of the operation $\widetilde{f}_r$ does produce rc-graphs (rather than just stable rc-graphs) because $R'$ has no crosses in positions $(r,w_r)$ and $(r+1,w_r)$. Hence $\widetilde{\partial}_r$ is a set-valued map from ${\mathcal R}(ws_r)$ to ${\mathcal R}(w)$.

\[
\begin{array}{c}
\mbox{\psfig{file=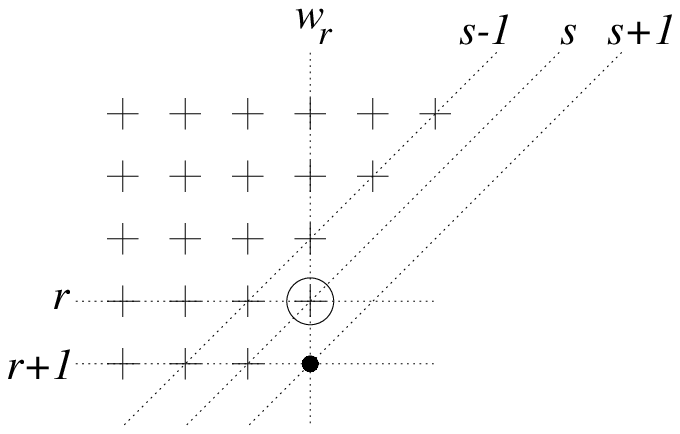}}
\end{array}
\]

Before stating the next result, we introduce some notation. Given positive integers $j>i$, let $(j..i)$ be the reduced word $j(j-1)\ldots i$. Given a word $u=u_1\ldots u_n$, we consider the following operator on words: $\phi_u:=u_1\sigma_1\ldots u_n\sigma_n$; in particular, $\phi(uv)=\phi_u(\phi(v))$, where $\phi$ is the plactification map (\ref{defpl}), and $uv$ is a reduced word with subwords $u$ and $v$. We need the following Remark. 

\begin{rema}\label{actphi}{\rm 
Let $v$ be a word with all letters strictly greater than $i$. The effect of applying $\phi_{(j..i)}$ to $v$ is to subtract 1 from the letters $i+1,\ldots,j+1$ in $v$, while leaving the other letters unchanged. Thus, the resulting word will not contain the letter $j+1$.}
\end{rema}

\begin{thm}\label{combk} We have that
\[{\mathcal R}(w)=\widetilde{\partial}_r({\mathcal R}(ws_r))\,,\;\;\;\;{\varPhi}(w)=\partial_r({\varPhi}(ws_r))\,,\]
where the right-hand sides are understood in the multiset sense discussed before Theorem \ref{combkey}.
\end{thm}
\begin{proof}
Any rc-graph $R'$ in ${\mathcal R}(w)$ has the form shown in the figure above, except that there is no cross in position $(r,w_r)$, and there might or might not be a cross in position $(r+1,w_r)$. If the former situation arises, this cross will clearly not be $r$-paired with a cross in row $r$. Hence it can be moved to row $r$ by successive applications of operations $\widetilde{e}_r$. Let $R:=\widetilde{e}_r^t(R')\cup\{(r,w_r)\}$, where $t$ is the number of unpaired crosses in row $r+1$ of $R'$. It is not hard to check that $R\in{\mathcal R}(ws_r)$, and that $R'\in\widetilde{\partial}_r(R)$. The uniqueness of $R$ with these properties is clear from the definition of $\widetilde{\partial}_r$.

Fix $R$ in ${\mathcal R}(ws_r)$, and let $s:=w_r+r-1$. In the above figure, the oblique dotted lines indicate the crosses in $R$ to which correspond identical letters (denoted $s-1$, $s$, $s+1$) in ${\rm red}(R)$. Let $u$, $u'$, $u''$, and $v$ be the subwords of ${\rm red}(R)$ corresponding to rows less than $r$, row $r$, row $r+1$, and rows greater than $r+1$, respectively. We have $u=\ldots(s-1..1)\ldots(s-1..r-1)$, $u'=\ldots(s..r)$, and $u''=\ldots\widehat{s+1}(s..r+1)$, where the $\widehat{\;\;}$ indicates the absence of a letter. Let $\widetilde{u}'$ be obtained from $u'$ by dropping the letter $s$. According to Remark \ref{actphi}, we have that $\phi_{(s-1..r)}(\phi(u''v))$ contains no letters $s$ and $s+1$. Therefore, this word is fixed by $\sigma_s$. Hence, $\phi(\widetilde{u}'u''v)$ is obtained from $\phi({u}'u''v)$ by dropping the letter $s$. By applying Remark \ref{actphi} repeatedly, we can see that the effect of applying $\phi_u$ on the letters $r,\ldots,s$ in the mentioned words is to subtract $r-1$ from these letters; furthermore, it is not hard to check that the effect on the other letters is identical in the two words. We conclude that $\phi(R\setminus\{(r,w_r)\})=\phi(R)\setminus\{(r,w_r)\}$. By applying Proposition \ref{commphi}, we can see that $\phi(\widetilde{\partial}_r(R))=\partial_r(\phi(R))$. 
\end{proof}

Hence, the set of diagrams ${\varPhi}(w)$ can be recursively generated via the operators $\partial_r$, starting from $\varPhi(w_0^n)$. This process is best illustrated using a tree, in a similar way to other recursive constructions in this paper. Let us note that different combinatorial versions of the divided difference operators were recently defined and used to recursively generate ${\mathcal R}(w)$ in \cite{milmrc}. On the other hand, one can easily show that $\varPhi(w_0^n)={\mathcal B}(w_0^n)=\{D(w_0^n)\}$; in fact, we can prove the more general statement in Proposition \ref{phibot} below. The next result follows by Schubert induction.

\begin{thm}\label{bij} Assuming Bergeron's result (\ref{combb}), we have that ${\mathcal B}(w)=\varPhi(w)$; hence, the map $\phi$ defines a bijection between ${\mathcal R}(w)$ and ${\mathcal B}(w)$. Assuming Winkel's result ${\mathcal K}(w)={\mathcal B}(w)$, we have that ${\mathcal K}(w)=\varPhi(w)$; hence, the map $\phi$ defines a bijection between ${\mathcal R}(w)$ and ${\mathcal K}(w)$. \end{thm}

The second statement was conjectured by Billey \cite{bilpri}. We intend to use the two characterizations of $\varPhi(w)$ presented above (the nonrecursive one in Remark \ref{impla} (1) and the recursive one in Theorem \ref{combk}) in order to find a transparent proof of Kohnert's conjecture, based on proving the stronger conjecture of Billey. The following two results are steps in this direction.

\begin{prop}\label{phibot}We have
\[\phi(R_{\rm bot}(w))=D(w)\,.\]
\end{prop}
\begin{proof}
Assuming that $w\in S_n$, we argue by induction on $n$. Let $R'$ be the rc-graph obtained from $R_{\rm bot}(w)$ by removing its first row, and let $\overline{w}:=w(R')$ be the corresponding permutation of $\{\,2,\ldots,n\,\}$. We have $w=w_1\widetilde{w}_2\ldots\widetilde{w}_n$, where
\[\case{\widetilde{w}_i}{:=}{\overline{w}_i-1}{2\le i\le w_1}{\overline{w}_i}\]
Thus, we can see that $D(w)$ is obtained from $D(\overline{w})$ by inserting an extra empty column $w_1$, as well as an extra first row containing $w_1-1$ left justified cells. By induction, we have that $\phi(R')=D(\overline{w})$. By Remark \ref{actphi}, the effect of the map $\phi_{{\rm red}_1(R)}=\phi_{(w_1-1..1)}$ on the word ${\rm col}(D(\overline{w}))$ is to subtract 1 from the letters $2,\ldots,w_1$. Therefore, the image is ${\rm col}(D({w}))$, which concludes the proof of the induction step.
\end{proof}

\begin{exa}\label{exbij}{\rm
Below are shown the Kohnert diagrams corresponding to the rc-graphs in Example \ref{exrc} via the map $\phi$. Note that the diagram on the left is just the diagram of the permutation $w=215463$. A sequence of Kohnert moves on this diagram which leads to the one on the right is: $4\rightarrow 2$, $2\rightarrow 1$, $5\rightarrow 4$, $4\rightarrow 2$. 
\[
\begin{array}{c}
\mbox{\psfig{file=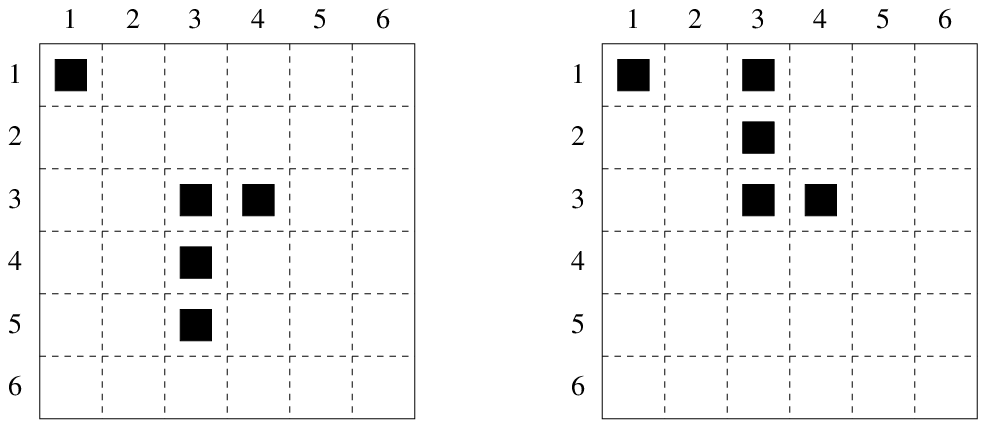}}
\end{array}
\]
}\end{exa}

We conclude this section by examining the special case of 321-avoiding permutations. Recall that a permutation $w$ is 321-avoiding if there are no $i<j<k$ such that $w_i>w_j>w_k$. Such permutations were first investigated in \cite{bjsscp}, where it was shown that the corresponding Schubert polynomials are {\em flagged skew Schur polynomials}; this fact is straightforwardly rederived in Remark \ref{balrem} (3) from the construction of a Schubert polynomial based on balanced labelings of the diagram of the corresponding permutation. Special cases of 321-avoiding permutations are Grassmannian permutations, which are permutations with a unique descent. The corresponding Schubert polynomials are Schur polynomials. 
If the permutation $w$ is 321-avoiding, the bijection $\phi$ can be made more explicit. Given a cross $(i,j)$ in an rc-graph, let us denote by ${\rm end}(i,j)$ the column where the line traversing it vertically ends. We need the following Lemma, which is best seen graphically, and is straightforward.

\begin{lem}\label{tto} If the permutation $w$ is 321-avoiding, then no rc-graph in ${\mathcal R}(w)$ has a line containing both a cross traversed vertically and one traversed horizontally. Conversely, if $w$ contains the pattern 321, then all rc-graphs in ${\mathcal R}(w)$ have such a line. 
\end{lem}

\begin{thm}\label{bij1} If $w$ is a 321-avoiding permutation, then the bijection $\phi$ maps every rc-graph $R$ to the diagram
\[\setp{(i,{\rm end}(i,j))}{(i,j)\in R}\,.\]
Furthermore, in this case we do have ${\mathcal K}(w)=\varPhi(w)$.
\end{thm}
\begin{proof}
We denote the above bijection from ${\mathcal R}(w)$ to some set of diagrams ${\mathcal E}(w)$ by $R\mapsto{\rm end}(R)$. Let us show first that ${\mathcal K}(w)={\mathcal E}(w)$. We start by noting that the order of the crosses in any row of an rc-graph $R$ coincides with the order of the corresponding cells in ${\rm end}(R)$; indeed, assuming that this is false contradicts the 321-avoiding property of $w$, by Lemma \ref{tto}. The inclusion ${\mathcal E}(w)\subseteq{\mathcal K}(w)$ is clear, because all L-moves on rc-graphs $R$ in ${\mathcal R}(w)$ move crosses up by just one row (by Lemma \ref{tto}), whence they become K-moves on ${\rm end}(R)$; the previous remark is also needed here, together with the obvious fact ${\rm end}(R_{{\rm bot}(w)})=D(w)$. The fact that ${\mathcal E}(w)$ is closed under K-moves is based on the following fact. If a cell is moved past $r$ cells in ${\rm end}(R)$ by a K-move, then the corresponding portion of $R$ looks as in the figure below (the remark at the beginning of the proof and Lemma \ref{tto} are needed again). Thus we can perform the $r+1$ indicated ladder moves in $R$, which together correspond to the considered K-move in ${\rm end}(R)$.
\[
\begin{array}{c}
\mbox{\psfig{file=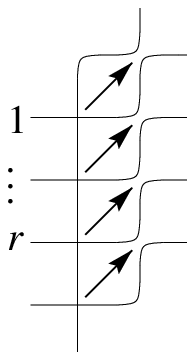}}
\end{array}
\] 

Assuming that $w\in S_n$, we now show that $\phi(R)={\rm end}(R)$ for any $R$ in ${\mathcal R}(w)$ by induction on $n$. Assume this is true for the rc-graph $R'$ obtained from $R_{\rm bot}(w)$ by removing its first row; the rows and columns of $R'$ are numbered $2,\ldots,n$. Write the word ${\rm red}_1(R)$ as $(j_p..i_p)\ldots(j_1..i_1)$ with $i_{k+1}-j_k\ge 2$ for $k=1,\ldots,p-1$. By Lemma \ref{tto}, the word ${\rm col}({\rm end}(R'))=\phi({\rm red}(R'))$ does not contain letters $i_1,\ldots,i_p$. Hence, by Remark \ref{actphi}, the effect of $\phi_{{\rm red_1(R)}}$ on this word is to subtract 1 from the letters $i_k+1,\ldots,j_k+1$, for $k=1,\ldots,p$, while leaving the other letters unchanged. We conclude that $\phi({\rm red}(R))={\rm col}({\rm end}(R))$. The example in the figure below illustrates our argument; here ${\rm red}_1(R)=(7..5)(2..1)$.
\[
\begin{array}{c}
\mbox{\psfig{file=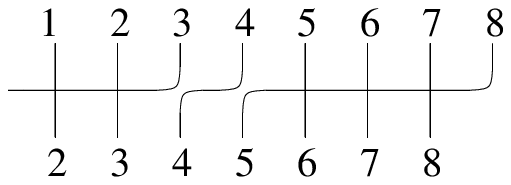}}
\end{array}
\]
\end{proof}

Note that, in general, the bijection $R\mapsto{\rm end}(R)$ is a monomial preserving bijection between ${\mathcal R}(w)$ and a different set of diagrams than ${\mathcal K}(w)$ (and $\varPhi(w)$).

\section{Balanced Labelings}

Fomin et al. \cite{fgrbls} present a construction of Schubert polynomials in terms of {\em balanced column-strict labelings} of the diagram of a permutation which satisfy a certain flag condition. A labeling of the cells of the diagram $D(w)$ is called balanced if it satisfies the following condition: if one rearranges the labels within any hook of $D(w)$ (which contains all the cells below and to the right of a cell of $D(w)$) such that they weakly increase from right to left and from top to bottom, then the corner label remains unchanged. A balanced labeling is called column-strict if no column contains two equal entries. The flag condition mentioned above is the requirement that all labels in row $i$ are at most $i$. For convenience, we will omit the terms ``column strict'' and ``flag condition'' for the rest of this section. With every balanced labeling is associated a monomial, obtained by taking the product of $x_i$ over all labels $i$. It was proved in \cite{fgrbls} that the Schubert polynomial $\fs_w(X)$ can be obtained by summing the monomials associated with all balanced labelings of $D(w)$.

In \cite{fgrbls} a monomial preserving bijection $R\mapsto B(R)$ was constructed from ${\mathcal R}(w)$ to balanced labelings of $D(w)$. The construction of the balanced labeling corresponding to an rc-graph $R$ for $w$  is based on the reduced word ${\rm red}(R)=a_1\ldots a_{l(w)}$ and the compatible sequence ${\rm comp}(R)=i_1\ldots i_{l(w)}$. The labeling is constructed as follows: the cell in row $p$ and column $q$ is labeled $i_r$ if the adjacent transposition $s_{a_r}=t_{a_r,a_r+1}$ in ${\rm red}(R)$ transposes $q$ and $w_p$, where $q<w_p$. Note that the number of labels $i$ is equal to the number of crosses in row $i$ of the rc-graph; this means that the bijection is monomial preserving. For instance, the balanced labeling of the diagram of $215463$  corresponding to the second rc-graph in Example \ref{exrc} is the following.
\[
\begin{array}{c}
\mbox{\psfig{file=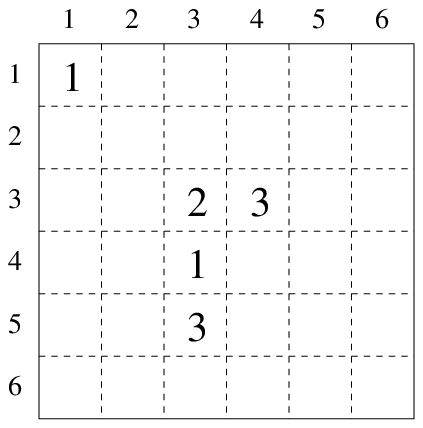}}
\end{array}
\]

We now intend to study the inverse map, from balanced labelings to rc-graphs. To this end, we associate a word $T(B)$ to a balanced labeling $B$, and a word $T(R)$ to an rc-graph $R$. The former is defined simply by sorting the entries in each row of $B$ in (weakly) decreasing order, and by concatenating the sequences obtained in this way, in increasing order of the rows. Now consider an rc-graph $R$ in ${\mathcal R}(w)$. Assume that the permutation $w=w_1\ldots w_n$ has $m$ descents, namely in positions $d_1,\ldots,d_m$; this means that $w_k>w_{k+1}$ if and only if $k=d_i$ for some $i$. We can associate a word $t_k=t_{k1}t_{k2}\ldots$ with line $k$ in $R$ in the following way: follow the line starting from its left end, and to each cross traversed horizontally associate a letter in $t_k$ equal to the number of the row where the cross lies. Now let $T_i:=t_{d_{i-1}+1}\ldots t_{d_i}$, where $i=1,\ldots,m$ and $d_0=0$; also let $T(R)=T_1\ldots T_m$. 

\begin{rem}\label{balrem}{\rm (1) The definition of the bijection $R\mapsto B(R)$ immediately implies that if $B=B(R)$, then $T(B)=T(R)$. Therefore, we can reconstruct $R$ from $T=T(B)$ only; that is, $R$ is determined by the {\em set} of labels in the rows of $B$. The precise description of the procedure is given in Algorithm \ref{bd2rc}.\\
\indent (2) If we sort the columns in a balanced labeling in ascending order, we obtain a sequence of words which can be read off from the corresponding rc-graph by following its lines starting at their top end, and by recording the row number of each cross traversed vertically. \\
\indent (3) The entries in the rows and columns of a balanced labeling are neither increasing nor decreasing, in general. In fact, based on Lemma \ref{tto} and Remark \ref{balrem} (1), we can easily check that if $w$ is 321-avoiding, then the entries in the rows of a balanced labeling of $D(w)$ are weakly decreasing, and the entries in columns are strictly increasing. It is well-known that, by omitting empty rows and columns in the diagram of a 321-avoiding permutation $w$, we obtain a skew Young diagram $\lambda/\mu=\lambda(w)/\mu(w)$ in French notation. Hence, the corresponding balanced labelings can be identified with the {\em reverse SSYT} of skew shape $\lambda/\mu$ (in English notation) satisfying a certain flag condition; the corresponding flag vector is clearly the decreasing sequence of integers $k$ for which the $k$th entry in ${\rm code}(w)$ in nonzero. Recall that a filling of a skew Young diagram is called a reverse SSYT (or column-strict plane partition) if the entries are weakly decreasing in rows and strictly decreasing in columns. In \cite{bjsscp}, these reverse SSYT are reinterpreted as SSYT of the shape obtained from $\lambda/\mu$ by a $180^\circ$ rotation. Thus, we have a straightforward way of rederiving the fact (proved in  \cite{bjsscp}) that Schubert polynomials indexed by 321-avoiding permutations are flagged skew Schur polynomials.\\
\indent (4) According to the previous Remark, if $w$ is a Grassmannian permutation with descent at $d$, then the map $R\mapsto B(R)$ followed by the renumbering $k\mapsto d+1-k$ of the entries in $B(R)$ is a bijection from ${\mathcal R}(w)$ to SSYT of shape $\lambda(w)$ with entries $1,\ldots,d$. This was already observed in \cite{berccs,bjsscp}. Note that the mentioned bijection is {\em not} monomial preserving; therefore, it is different from the map $R\mapsto \widetilde{Q}(R)$ (discussed in Section 4) between the same sets as above, which is a bijection for Grassmannian permutations. Hence, it is natural to ask about the way in which the SSYT $\widetilde{Q}(R)$ is related to the reverse SSYT $B(R)$ of the same shape. This is addressed in Theorem \ref{revtab} below.
}
\end{rem}

\begin{thm}\label{revtab} For any Grassmannian permutation $w$ and $R$ in ${\cal R}(w)$, we have
\[\widetilde{Q}(R)={\rm evac}(B(R))\,.\]
More generally, if $w$ is 321-avoiding, then $\widetilde{Q}(R)$ is the rectification of the SSYT of skew shape obtained from $B(R)$ by a $180^\circ$ rotation. 
\end{thm}
\begin{proof}
Let $w$ be 321-avoiding. By Theorem \ref{pla} (3) and (\ref{symm}), we have
\[\widetilde{Q}(R)=Q(\phi(R))=P(\phi(R)')\,;\]
recall that 
\[\phi(R)'=\binom{\phi({\rm red}(R))'}{{\rm comp}(R)'}\,,\]
cf. (\ref{defphi}) and the definition of the biword $w'$ (corresponding to some biword $w$) at the beginning of Section 3. On the other hand, we know from Theorem \ref{bij1} that $\phi(R)={\rm end}(R)$. Therefore, by Remark \ref{balrem} (2), we have that $({\rm comp}(R)')^{\rm rev}$ is precisely the column word of $B(R)$. In conclusion, $\widetilde{Q}(R)$ can be obtained by applying Schensted insertion to the reverse column word of $B(R)$, and therefore coincides with ${\rm evac}(B(R))$ when $w$ is Grassmannian.
\end{proof}

\begin{alg}\label{bd2rc}$\;\;$\newline
\noindent for $k=1,\ldots,n$ do\newline
\indent start at the left end $(k,1)$ of line $k$;\newline
\indent while not at the upper end of line $k$ do\newline
\indent\indent let $(i,j)$ be the current position;\newline
\indent\indent if a previous line crosses position $(i,j)$ horizontally then\newline
\indent\indent\indent cross position $(i,j)$ vertically;\newline
\indent\indent else\newline
\indent\indent\indent if $(i,j)\not =(k,1)$ then avoid position $(i,j)$ fi;\newline
\indent\indent\indent move to the right past a number of lattice points\newline
\indent\indent\indent\indent\indent equal to the number of labels $i$ in the $k$th row of $T$;\newline
\indent\indent\indent avoid the current position;\newline
\indent\indent fi;\newline
\indent od;\newline
od.
\end{alg}

We conclude with a closer look at the word $T(R)$ associated to an rc-graph $R$. We view $T(R)$ as the sequence of its subwords $T_i=t_{d_{i-1}+1}\ldots t_{d_i}$, where $i=1,\ldots,m$. Note that the length of the subword $t_k$ is the $k$th entry $c_k$ in the code of $w(R)$. Hence the lengths of $t_{d_{i-1}+1},\ldots, t_{d_i}$ form a weakly increasing sequence. Now view $T_i$ as a filling of a Young diagram of shape $(c_{d_{i-1}+1},\ldots, c_{d_i})$, where the French convention for such diagrams is used. 

\begin{prop} Using the notation above, we have that $T(R)$ is a sequence of reverse SSYT \linebreak $(T_1,\ldots,T_m)$.\end{prop}
\begin{proof}
Fix $i$, $k$ such that $d_{i-1}+1\le k<k+1\le d_i$, and consider the lines $k$ and $k+1$ in $R$. It is not hard to see that all lines which cross line $k$ vertically will have to cross line $k+1$ vertically as well; indeed, such lines cannot cross line $k$ twice, so they must exit the region between lines $k$ and $k+1$ by crossing the latter. Now remove all crosses in the mentioned region. By following the lines which cross line $k$ vertically in the new diagram, we obtain a pairing of these crosses with some crosses on line $k+1$ (see the figure below). Obviously, the row of the first cross in any pair is strictly smaller than the row of the second. Furthermore, we can order the pairs such that the first and second components appear in the order obtained by following lines $k$ and $k+1$ from their left ends to their top ends. This essentially shows that the entries in the columns of $T_i$ strictly increase.   
\[
\begin{array}{c}
\mbox{\psfig{file=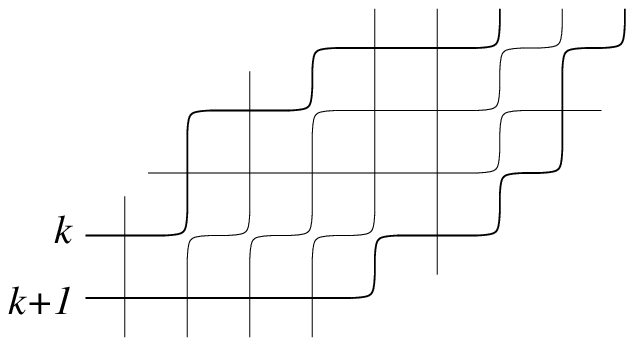}}
\end{array}
\]
\end{proof}

\begin{exa}{\rm
The sequences of reverse SSYT corresponding to the two rc-graphs in Example \ref{exrc} are shown below, stacked on top of each other. More precisely, there are three reverse SSYT, corresponding to the following sets of rows: $\{\,1\,\}$, $\{\,2,3\,\}$, $\{\,4,5\,\}\,$.
\[\begin{array}{cccccc}
1&$\;$&$\;\;\;\;\;\;\;\;\;\;\;\;\;\;$&1$\;$&$\;$&$\;$\\
$\,$\\
3&3&$\;$&3&2\\
4&$\;$&$\;$&1&$\;$\\
5&$\;$&$\;$&3&$\;$\\\end{array}\]
}\end{exa}

\begin{rema}{\rm The sequence of reverse SSYT $T(R)=(T_1,\ldots,T_m)$ associated to an rc-graph $R$ might be useful for the important open problem of multiplying Schubert polynomials. One approach to this problem, adopted in \cite{babrcg,koggsi,kakppf}, is to define an insertion procedure for rc-graphs (i.e., to ``insert'' an rc-graph into another) which generalizes Schensted insertion. There are such procedures (see the mentioned papers) which ``insert'' rc-graphs $R$ corresponding to certain Grassmannian permutations. They can be viewed as ``insertions'' of the corresponding words $T(R)$, and they depend on a parameter $d$, which is the unique descent of the permutation $w(R)$; we call such an insertion a $d$-insertion. Based on the results of the mentioned papers, we suggest that the general procedure of ``inserting'' an arbitrary rc-graph $R$ into another could be reduced to a sequence of $d_i$-insertions of the words $T_i$ corresponding to $R$.
}\end{rema}


\newpage

\[
\begin{array}{c}
\mbox{\psfig{file=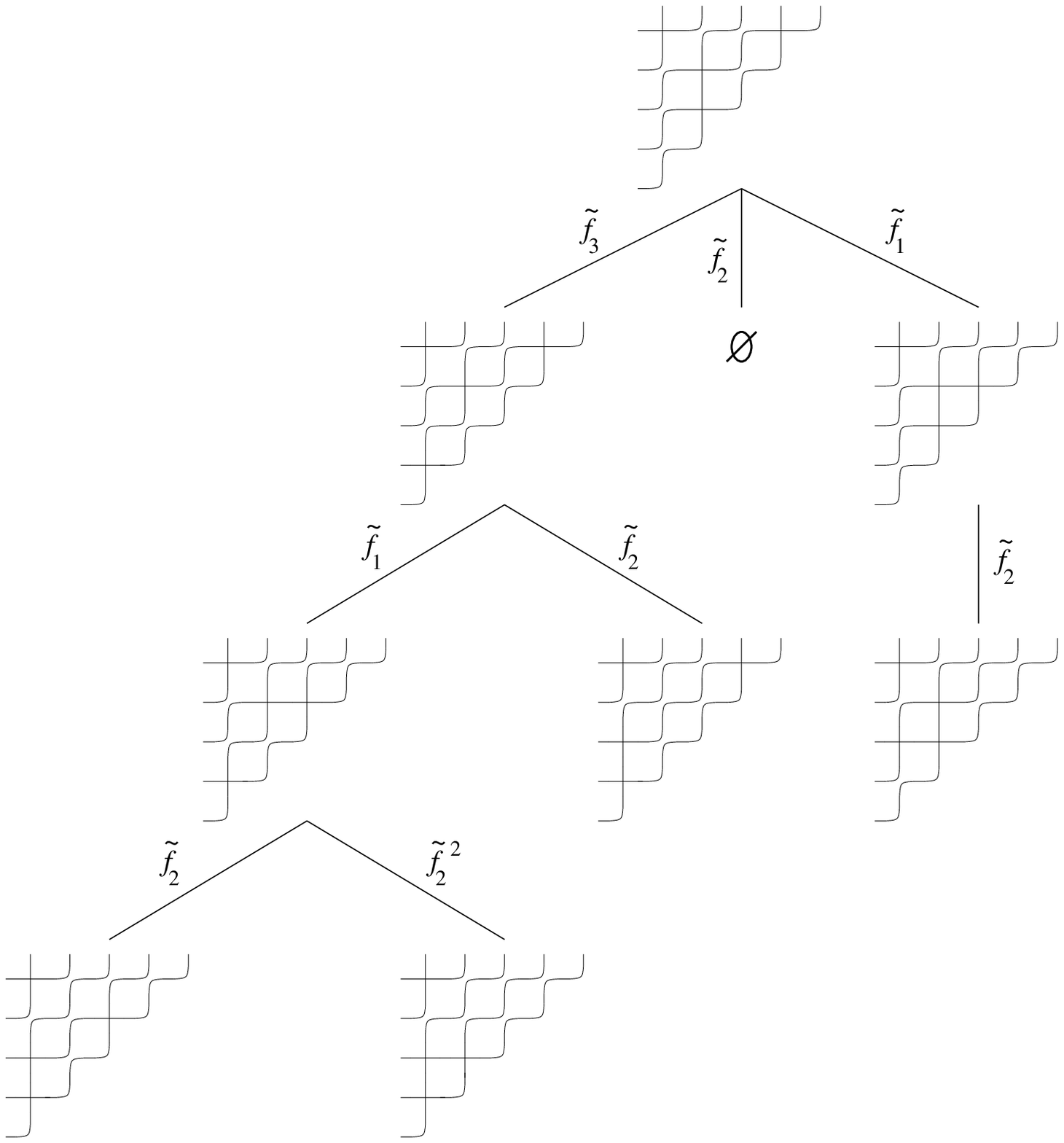}}
\end{array}
\]

Appendix. The construction of the crystal in Example \ref{excrys}.

\vspace{5mm}

\end{document}